\documentclass[12pt]{article}

\usepackage{fontspec, xunicode, xltxtra}  
\usepackage[slantfont,boldfont]{xeCJK}
\usepackage{mathtools}
\usepackage{amsmath,diagbox}
\usepackage{indentfirst}
\usepackage{mathrsfs}
\usepackage{amsfonts}
\usepackage{arydshln}
\usepackage{enumerate}
\usepackage{setspace}
\usepackage{amssymb,amsthm,cases}
\usepackage{amsbsy,pifont}
\usepackage{latexsym,diagbox,tabularx}
\usepackage{amsmath,amscd,bbm,multirow}
\usepackage{graphicx,epsfig,extarrows,dsfont,mathtools}
\usepackage[final,allcolors=blue,colorlinks=true]{hyperref}
\usepackage[normalem]{ulem}
\input xy
\usepackage{caption,subcaption} 
\usepackage{tikz}
\usepackage{tikz-cd}
\usetikzlibrary{decorations.pathreplacing,decorations.markings}
\usetikzlibrary{cd}
\usetikzlibrary{shapes.geometric}
\usetikzlibrary{arrows,arrows.meta}
\usetikzlibrary{graphs,positioning}
\usetikzlibrary{matrix,decorations.pathmorphing}
\usetikzlibrary{math,calc,intersections,through,angles,arrows.meta,shapes.geometric,shadows,quotes,spy,datavisualization,datavisualization.formats.functions,plotmarks}
\tikzset{every picture/.style={samples=300,smooth,line join=round,thick,>=stealth}}
\usepackage{enumerate}
\usepackage{fancyhdr}
\usepackage[sectionbib]{chapterbib}

\tikzset{
	on each segment/.style={
		decorate,
		decoration={
			show path construction,
			moveto code={},
			lineto code={
				\path [#1]
				(\tikzinputsegmentfirst) -- (\tikzinputsegmentlast);
			},
			curveto code={
				\path [#1] (\tikzinputsegmentfirst)
				.. controls
				(\tikzinputsegmentsupporta) and (\tikzinputsegmentsupportb)
				..
				(\tikzinputsegmentlast);
			},
			closepath code={
				\path [#1]
				(\tikzinputsegmentfirst) -- (\tikzinputsegmentlast);
			},
		},
	},
	mid arrow/.style={postaction={decorate,decoration={
				markings,
				mark=at position .5 with {\arrow[#1]{stealth}}
	}}},
}

\numberwithin{equation}{section}

\theoremstyle{plain}
\newtheorem{theorem}{Theorem}[section]

\newtheorem{lemma}[theorem]{Lemma}
\newtheorem{proposition}[theorem]{Proposition}

\theoremstyle{definition}

\newtheorem{remark}[theorem]{Remark}

\def\XXint#1#2#3{{\setbox0=\hbox{$#1{#2#3}{\int}$}
		\vcenter{\hbox{$#2#3$}}\kern-.5\wd0}}

\DeclareMathSymbol{\subseteq}{\mathrel}{symbols}{"12}
\DeclareMathSymbol{\supseteq}{\mathrel}{symbols}{"13} 
\DeclareMathSymbol{\subsetneq}{\mathrel}{AMSb}{"28}                  \DeclareMathSymbol{\supsetneq}{\mathrel}{AMSb}{"29}    
 
\DeclareMathSymbol{\nsubseteq}{\mathrel}{AMSb}{"2A}                  

\DeclareMathSymbol{\nsupseteq}{\mathrel}{AMSb}{"2B}

\DeclareMathDelimiter{\langle}{\mathop}{symbols}{"68}{largesymbols}{"0A}
\DeclareMathDelimiter{\rangle}{\mathclose}{symbols}{"69}{largesymbols}{"0B}
\DeclareSymbolFont{txfontsA}{U}{txmia}{m}{it}
\SetSymbolFont{txfontsA}{bold}{U}{txmia}{bx}{it}
\DeclareFontSubstitution{U}{txmia}{m}{it}
\DeclareMathSymbol{\upalpha}{\mathord}{txfontsA}{"0B}
\DeclareMathSymbol{\upbeta}{\mathord}{txfontsA}{"0C}
\DeclareMathSymbol{\upgamma}{\mathord}{txfontsA}{"0D}
\DeclareMathSymbol{\updelta}{\mathord}{txfontsA}{"0E}
\DeclareMathSymbol{\upepsilon}{\mathord}{txfontsA}{"0F}
\DeclareMathSymbol{\upzeta}{\mathord}{txfontsA}{"10}
\DeclareMathSymbol{\upeta}{\mathord}{txfontsA}{"11}
\DeclareMathSymbol{\uptheta}{\mathord}{txfontsA}{"12}
\DeclareMathSymbol{\upiota}{\mathord}{txfontsA}{"13}
\DeclareMathSymbol{\upkappa}{\mathord}{txfontsA}{"14}
\DeclareMathSymbol{\uplambda}{\mathord}{txfontsA}{"15}
\DeclareMathSymbol{\upmu}{\mathord}{txfontsA}{"16}
\DeclareMathSymbol{\upnu}{\mathord}{txfontsA}{"17}
\DeclareMathSymbol{\upxi}{\mathord}{txfontsA}{"18}
\DeclareMathSymbol{\uppi}{\mathord}{txfontsA}{"19}
\DeclareMathSymbol{\uprho}{\mathord}{txfontsA}{"1A}
\DeclareMathSymbol{\upsigma}{\mathord}{txfontsA}{"1B}
\DeclareMathSymbol{\uptau}{\mathord}{txfontsA}{"1C}
\DeclareMathSymbol{\upupsilon}{\mathord}{txfontsA}{"1D}
\DeclareMathSymbol{\upphi}{\mathord}{txfontsA}{"1E}
\DeclareMathSymbol{\upchi}{\mathord}{txfontsA}{"1F}
\DeclareMathSymbol{\uppsi}{\mathord}{txfontsA}{"20}
\DeclareMathSymbol{\upomega}{\mathord}{txfontsA}{"21}
\DeclareMathSymbol{\upvarepsilon}{\mathord}{txfontsA}{"22}
\DeclareMathSymbol{\upvartheta}{\mathord}{txfontsA}{"23}
\DeclareMathSymbol{\upvarpi}{\mathord}{txfontsA}{"24}
\DeclareMathSymbol{\upvarrho}{\mathord}{txfontsA}{"25}
\DeclareMathSymbol{\upvarsigma}{\mathord}{txfontsA}{"26}
\DeclareMathSymbol{\upvarphi}{\mathord}{txfontsA}{"27}                   

\DeclareSymbolFont{ugmL}{OMX}{mdugm}{m}{n}
\SetSymbolFont{ugmL}{bold}{OMX}{mdugm}{b}{n}
\DeclareMathAccent{\wideparen}{\mathord}{ugmL}{"F3}

\def\pt{\partial}

\def\ra{\rightarrow}

\def\s{\subseteq}

\def\ol{\overline}

\def\vp{\varphi}

\def\es{\emptyset}

\def\pt{\partial}

\def\Om{\Omega}
\def\la{\lambda}
\def\al{\alpha}
\def\be{\beta}
\def\de{\delta}
\def\ga{\gamma}

\def\si{\sigma}

\def\La{\Lambda}

\def\ts{\times}
\def\ka{\kappa}

\def\iy{\infty}

\def\f{\frac}

\def\vs{\vspace}

\def\df{\mathrm d}
\def\im{\mathrm i}

\def\wh{\widehat}

\def\hra{\hookrightarrow}

\def\mcL{\mathcal{L}}

\def\mcB{{\mathcal{B}}}

\def\mcA{\mathcal{A}}

	\def\Ee{\mathbbm{1}}

	\def\mbT{\mathbb{T}}

	\DeclareMathOperator{\sgn}{sgn}

	\DeclareMathOperator{\Real}{Re}
	\DeclareMathOperator{\Image}{Im}

	 \DeclareMathOperator{\Ran}{Ran}

	\newcommand{\R}{\mathbb R}
	\newcommand{\N}{\mathbb N}
	\newcommand{\Z}{\mathbb Z}
	
	\newcommand{\Cc}{\mathbb C}

\begin{document} 
	
	\title{Boundary Observability for a Class of Degenerate Kolmogorov-Type Equations in a Cylindrical Domain}

	\author{ Dong-Hui Yang$^{a}$, Mengze Gu$^{b}$ \footnote{\small
			The corresponding author: gu3059268737@gmail.com}, Yuanzhi Zhou$^{a}$,
		Ghadir Shokor$^{c}$
		\\
		$^a${\it School of Mathematics and Statistics, Central South University}\\
		{\it Changsha 410075, P.R.China}\\
		$^b${\it School of Mathematics and Science, Changsha Normal University}\\ {\it  Changsha 410083, China}\\
		$^c${\it School of Mathematics and Statistics, Beijing Institute of Technology}\\
		{\it Beijing 100081, P.R.China} }

	\date{}
	\maketitle{}
	\thispagestyle{empty}
	\thispagestyle{empty}
	
	\begin{abstract}
This paper investigates the boundary observability for a class of multi-dimensional degenerate Kolmogorov-type kinetic equations in a cylindrical domain $\Omega = \mathbb{T}_{x_1} \times (0,1)_{x_2}$. The model features a horizontal drift term $x_2^\gamma \partial_{x_1} \varphi$ coupled with a boundary-degenerate diffusion operator $\partial_{x_2}(x_2^\alpha \partial_{x_2} \varphi)$ acting solely on the velocity variable, where $\gamma = 1$ and $\alpha \in (0,1)$. By leveraging the structural advantages of the cylindrical domain, we establish an analytical framework based on the method of separation of variables. Crucial Carleman estimates tailored to the degenerate dual operator are derived to capture the subtle dissipation effects under the influence of the transport term. Consequently, by virtue of a uniform a priori estimate, we establish the quantitative boundary observability inequality for the state function. 
		
\vspace{3mm} \noindent{\textbf{Keywords:} Boundary observability, Degenerate Kolmogorov-type equation, Separation of variables, Carleman estimates, Cylindrical domain.} 
		
\vs{3mm}
		
		
\end{abstract}

\section{Introduction}

The controllability of partial differential equations (PDEs) has been extensively investigated in a large body of literature (see, e.g., \cite{Coron,Fursikov,Lions,Lu,Privat,Tucsnak,Yao,Zuazua}). 

The study of controllability for degenerate partial differential equations has also attracted significant interest \cite{Akil,Alabau,Araruna,Bai,Beauchard1,Beauchard,Beauchard3,Beauchard2,Buffe,Cannarsa1,Fragnelli,Fragnelli1,Gueye,Guo,Martinez,Zhang}. However, most of these works focus exclusively on the one-dimensional case, with notable exceptions such as those involving the Grushin operator \cite{Beauchard3}. Recently, several results have been established in higher dimensions ($N\geq 2$) \cite{Yang5,Yang,Yang4,Yang6,Yang3,Yang8,Yang1,Yang7}, employing various techniques such as the cut-off method, approximation methods, the shape design method, and the separation of variables. 

According to J.-L. Lions' Hilbert Uniqueness Method (HUM), the controllability of a linear PDE is equivalent to the observability of its adjoint (dual) system \cite{Lions}. This fundamental equivalence also extends to degenerate partial differential equations \cite{Yang2}. 

In this paper, we focus on the boundary observability for a class of degenerate Kolmogorov-type equations. 

Let $\mbT=\R/(2\pi\Z)$ be the one-dimensional torus, represented by the angular coordinate $x_1\in[0,2\pi)$ with periodic identification and measure $\mathrm{d}x_1$. Set $\La=(0,1)$,  $\Om=\mbT_{x_1}\times\La_{x_2}$, $Q=\Omega\ts (0,T)$, and $\Sigma=\pt\Om\ts (0,T)$ with $\pt\Om=\mbT_{x_1}\ts \{0,1\}$. For $\ga=1$ and $\alpha\in(0,1)$, we consider the following degenerate kinetic equation:
\begin{equation}\label{10.31.1}
	\begin{cases}
		\partial_t \vp - x_2^\ga\partial_{x_1} \vp - \partial_{x_2}\left(x_2^\alpha\partial_{x_2} \vp\right) = u, &\mbox{in } Q,\\[1mm]
		\vp = 0, &\mbox{on }\Sigma,\\[1mm]
		\vp(0) = \vp^0, &\mbox{in }\Om,
	\end{cases}
\end{equation}
where $\vp^0\in L^2(\Om)$ is the initial datum and $u\in L^2(Q)$ acts as the control function. Let us denote the associated spatial operator by
\begin{equation}\label{10.31.7}
	\mcA \vp = -\left(x_2^\ga \pt_{x_1}\vp + \pt_{x_2}\left(x_2^\al \pt_{x_2}\vp\right)\right). 
\end{equation}

Equation \eqref{10.31.1} belongs to the class of degenerate kinetic equations (or Kolmogorov-type equations). It features a horizontal transport term $x_2^\gamma\partial_{x_1} \varphi$ coupled with a boundary-degenerate diffusion operator $\partial_{x_2}(x_2^\alpha\partial_{x_2} \varphi)$ acting solely on the velocity variable $x_2$. 

The controllability and observability problems for various Kolmogorov equations have been extensively investigated in the literature (see, e.g., \cite{Beauchard1,Beauchard,Beauchard2} and the references therein). However, the specific degenerate principal operator in \eqref{10.31.1} in higher dimensions ($N \ge 2$) has not been previously addressed. Inspired by the pioneering approaches in \cite{Beauchard,Beauchard2}, this paper aims to fill this gap. To overcome the additional analytical difficulties arising from the spatial degeneracy, existing techniques—most notably the Carleman estimates—must be substantially refined; accordingly, we develop and employ a different parameter-weighted Carleman estimate tailored to our degenerate setting.

In this work, we shall primarily focus on the case $\gamma = 1$. The reasons for this choice are two-fold:

First, the case $\gamma=1$ gives a linear shear transport term $x_2\partial_{x_1}$. The coupling of horizontal transport and vertical diffusion is also present in linearized Crocco-type models; see \cite{Martinez}. Here we study this coupling with the specific boundary-degenerate diffusion coefficient $x_2^\alpha$ and periodic horizontal geometry.

Second, establishing boundary observability in higher dimensions ($N \ge 2$) already poses substantial mathematical challenges due to the subtle interplay between the boundary degeneracy $x_2^\alpha$ and the drift term. Focusing on $\gamma = 1$ allows us to present our core methodology---such as the separation of variables and spectral decoupling---in a clear and transparent manner, avoiding excessive algebraic technicalities.

The modal enhanced dissipation estimate is established below for $\gamma\in\{1,2\}$, and the Carleman estimate is stated for fixed $\gamma\in\mathbb N^*$. The global boundary observability theorem is proved for $\gamma=1$. The limitation of the present proof when $\gamma=2$ is discussed in Remark~\ref{11.07.R2}.

The rest of this paper is organized as follows. 
In Section \ref{S2}, we introduce the relevant weighted Sobolev spaces, present the well-posedness (existence and uniqueness) of the mild solution to equation \eqref{10.31.1}, and describe the analytical framework based on the method of separation of variables. 
Section \ref{S3} is devoted to the derivation of the crucial Carleman estimates for the dual system. 
In Section \ref{S4}, by virtue of a uniform a priori estimate, we establish the quantitative boundary observability inequality. 

The main mathematical contributions of this work are crystallized in Lemmas \ref{10.31.L3}, \ref{10.31.L4}, and \ref{11.05.L2}, Proposition \ref{11.04.P1}, and ultimately, the main observability result formulated in Theorem \ref{11.07.T2}.

\section{Preliminaries}\label{S2}

In this section, we introduce the weighted Sobolev spaces tailored to the boundary degeneracy, outline the framework for the method of separation of variables, and establish the well-posedness (existence and uniqueness) of the mild solution to the degenerate system \eqref{10.31.1}.

\subsection{Weighted Sobolev spaces}
Define
\begin{equation*}
	V=\left\{\psi\in L^2(\Om)\colon   x_2^\f{\al}{2} \pt_{x_2}\psi\in L^2(\Om)\right\}.
\end{equation*}
The inner product and norm on $V$ are defined by
\begin{equation*}
	\begin{split}
		(\vp,\psi)_V=\int_\Om \Big(\vp\psi+x_2^\al (\pt_{x_2} \vp)(\pt_{x_2}\psi)\Big)\df x, \mbox{ and } \|\psi\|_V=(\psi,\psi)_V^\f{1}{2}. 
	\end{split}
\end{equation*}

Set 
\begin{equation*}
	\begin{split}
		V_0=\mbox{the closure of $C_0^\iy(\Om)$ in } V, \mbox{ and } V_0'=\mbox{the dual space of }V_0. 
	\end{split}
\end{equation*}

According to \cite{GC, Heinonen}, the spaces $V$, $V_0$, and $V_0'$ are well-defined separable Hilbert spaces.

Define 
\begin{equation*}
	D(\mcA)=\left\{\psi\in V_0\colon \mcA\psi\in L^2(\Om)\right\}. 
\end{equation*}

\begin{lemma}[{\cite[Lemma 2.1]{Yang5}}]\label{10.31.L1}
	Let $\al\in (0,1)$. For each $\psi\in V_0$, we have 
	\begin{equation}\label{10.31.2}
		\int_\Om x_2^{\al-2}\psi^2\df x\leq \f{4}{(1-\al)^2}\int_\Om x_2^\al (\pt_{x_2}\psi)^2\df x.
	\end{equation} 
\end{lemma}

\begin{proof}
	By the density argument, we only need to show the case $\psi\in C_0^\iy(\Om)$. 
	
	For each $\be\in (\al,1)$, from 
	\begin{equation*}
		\begin{split}
			\int_s^1 x_2^{\al-1-\be}\df x_2=\left.\f{1}{\al-\be}x_2^{\al-\be}\right|_{x_2=s}^{x_2=1}\leq \f{1}{\be-\al}s^{\al-\be}, 
		\end{split}
	\end{equation*}
	we get
	\begin{equation*}
		\begin{split}
			\int_\Om x_2^{\al-2}\psi^2\df x
			&=\int_\mbT \int_0^1 x_2^{\al-2} (\psi(x_1,x_2)-\psi(x_1,0))^2\df x_2\df x_1\\
			&=\int_\mbT \int_0^1x_2^{\al-2}\left(\int_0^{x_2} \f{\pt \psi}{\pt x_2}(x_1,s)\df s\right)^2\df x_2\df x_1\\
			&\leq \int_\mbT\int_0^1 x_2^{\al-2}\left(\int_0^{x_2} s^{\be} \left|\f{\pt \psi}{\pt x_2}\right|^2(x_1,s)\df s\right)\left(\int_0^{x_2}s^{-\be}\df s\right)\df x_2\df x_1 \\
			&\leq \f{1}{1-\be}\int_\mbT \int_0^1 \int_0^{x_2} x_2^{\al-1-\be} s^\be \left|\f{\pt\psi}{\pt x_2}\right|^2(x_1,s)\df s\df x_2\df x_1\\
			&=\f{1}{1-\be}\int_\mbT \int_0^1\int_s^1x_2^{\al-1-\be}s^\be \left|\f{\pt \psi}{\pt x_2}\right|^2(x_1,s)\df x_2\df s\df x_1\\
			&\leq \f{1}{(\be-\al)(1-\be)}\int_\mbT \int_0^1 s^\al \left|\f{\pt \psi}{\pt x_2}\right|^2(x_1,s)\df s\df x_1. 
		\end{split}
	\end{equation*}
	Taking $\be=\f{1+\al}{2}$, we obtain \eqref{10.31.2}.  
\end{proof}

\begin{remark}\label{10.31.R1}
	From the inequality \eqref{10.31.2}, we obtain 
	\begin{equation}\label{10.31.3}
		\int_\Om \psi^2\df x=\int_\Om x_2^{2-\al}\left(x_2^{\al-2}\psi^2\right)\df x\leq \int_\Om x_2^{\al-2}\psi^2\df x\leq \f{4}{(1-\al)^2}\int_\Om x_2^\al (\pt_{x_2}\psi)^2\df x. 
	\end{equation}

This is a weighted Poincar\'e inequality on $V_0$. On $V_0$ we use the equivalent energy norm
\begin{equation*}
\|\psi\|_{V_0}	:=\left(\int_\Omega x_2^\alpha |\partial_{x_2}\psi|^2\,\mathrm{d}x\right)^{1/2}.
\end{equation*}
This norm is equivalent to the restriction of $\|\cdot\|_V$ to $V_0$.

It is worth emphasizing that the positive embedding constant $\frac{4}{(1-\alpha)^2}$ is completely independent of the transport exponent $\ga$. 
	
Note that due to the geometry of the domain and the degenerate nature of the weight, the continuous embedding $V_0 \hookrightarrow L^2(\Om)$ is not compact, although we have the standard chain of inclusions $H_0^1(\Om) \s  V_0 \s  L^2(\Om)$.
\end{remark}

\subsection{Method of separation of variables: step 1}\label{S2.2}

Define  
\begin{equation*}
	H^1(\La;w)=\left\{u\in L^2(\La)\colon x_2^\f{\al}{2} \pt_{x_2}u\in L^2(\La) \right\}, \mbox{ with } w=x_2^\al. 
\end{equation*}
Its inner product and norm are defined by 
\begin{equation*}
	(u,v)_{H^1(\La;w)}=\int_\La \left(uv+w(\pt_{x_2}u)(\pt_{x_2}v)\right)\df x_2, \mbox{ and } \|u\|_{H^1(\La;w)}=(u,u)_{H^1(\La;w)}^\f{1}{2}. 
\end{equation*}
Set 
\begin{equation*}
	\begin{split} 
	H_0^1(\La;w)
	&=\mbox{the closure of $C_0^\iy(\La)$ in }H^1(\La;w),\\ H^{-1}(\La;w)
	&=\mbox{the dual space of $H_0^1(\La;w)$ with pivot space } L^2(\La). 
	\end{split} 
\end{equation*}

According to \cite{GC, Heinonen}, the spaces $H_0^1(\La;w), H^1(\La;w)$ and $H^{-1}(\La;w)$ are well-defined separable Hilbert spaces. Moreover, from Lemma \ref{10.31.L1} and Remark \ref{10.31.R1}, we obtain 

\begin{remark}\label{10.30.R2}
	(1) (see \cite[(2.1), p.~165]{Alabau}, or \cite[Lemma 2.1]{Yang}.) For each $u\in H_0^1(\La;w)$, we have 
	\begin{equation*}
		\int_\La x_2^{\al-2}u^2\df x_2\leq \f{4}{(1-\al)^2}\int_\La x_2^\al (\pt_{x_2}u)^2\df x_2. 
	\end{equation*}
	Moreover, 
	\begin{equation*}
		\int_\La u^2\df x_2 \leq \f{4}{(1-\al)^2}\int_\La x_2^\al (\pt_{x_2}u)^2\df x_2. 
	\end{equation*}
	
(2) From now on, we equip $H_0^1(\La;w)$ with the equivalent energy norm
\begin{equation*}
\|u\|_{H_0^1(\La;w)}	:=\left(\int_\La x_2^\alpha	|\partial_{x_2}u|^2\,\mathrm{d}x_2\right)^{1/2}.
\end{equation*}

	(3) The embedding $H_0^1(\La;w)\hra L^2(\La)$ is compact (see  \cite[Corollary 3.1 (ii)]{Yang6}). 
\end{remark}

From this point onward we use the complexifications of the
weighted Hilbert spaces introduced above. Inner products are
linear in the first argument and conjugate-linear in the second;
in particular,
\[
(u,v)_{L^2(\Omega)}=\int_{\Om} u\overline v \df x,\qquad 
(u,v)_V=\int_\Omega\left(u\overline v+
x_2^\alpha\partial_{x_2}u\,
\overline{\partial_{x_2}v}\right)\mathrm{d}x.
\]
All energy identities below are understood after taking real parts.
The real-valued solutions are recovered by restricting the
complexified evolution to real initial data and forcing.

To perform a rigorous spectral analysis of $\mcA$, we leverage the geometric structure of the domain $\Om = \mbT_{x_1} \times \La$. Since the one-dimensional torus $\mbT_{x_1}$ is compact, any function $\vp \in L^2(\Om)$ can be decomposed into a Fourier series with respect to the spatial variable $x_1$:
\begin{equation}\label{10.31.11}
	\vp(x_1, x_2) = \sum_{k \in \mathbb{Z}} \widehat{\vp}_k(x_2) \left(\f{1}{\sqrt{2\pi}}e^{\im k  x_1}\right),
\end{equation}
where the Fourier coefficients $\widehat{\vp}_k(x_2)$ belong to $L^2(\La)$ for each $k \in \mathbb{Z}$. This structural decomposition is justified by the tensor product property of the underlying Hilbert spaces. Specifically, since the domain $\Omega = \mbT_{x_1} \times \La$ is a product space, $L^2(\Omega)$ can be isometrically identified with the Hilbert tensor product $L^2(\mbT_{x_1}) \otimes L^2(\La)$. Recalling that the family of complex exponentials $\{\f{1}{\sqrt{2\pi}}e^{\im k  x_1}\}_{k \in \mathbb{Z}}$ forms a complete orthonormal basis of $L^2(\mbT_{x_1})$, any function $\vp \in V_0 \s  L^2(\Omega)$ uniquely admits the expansion \eqref{10.31.11}. By Fubini's theorem and Parseval's identity, the Fourier coefficients defined by
\begin{equation}\label{10.31.14}
	\widehat{\vp}_k(x_2) = \frac{1}{\sqrt{2\pi}} \int_{\mbT_{x_1}} \vp(x_1, x_2) e^{-\im k  x_1} \,\mathrm{d}x_1
\end{equation}
satisfy $\sum_{k \in \mathbb{Z}} \|\widehat{\vp}_k\|_{L^2(\La)}^2 = \|\vp\|_{L^2(\Omega)}^2 < +\infty$, ensuring the strong convergence of the series in $L^2(\Omega)$.

The Fourier decomposition \eqref{10.31.11} allows us to decouple the global operator $\mcA$ into a family of one-dimensional localized operators $\{\mcA_k\}_{k \in \mathbb{Z}}$ acting solely on the velocity variable $x_2$. That is, for each $k \in \mathbb{Z}$, the operator $\mcA_k$ is defined by
\begin{equation}\label{10.31.12}
	\mcA_k \colon D(\mcA_k) \s  L^2(\La) \to L^2(\La),\quad \mcA_k \psi = -\im k x_2^\ga \psi - \partial_{x_2}\left(x_2^\alpha \partial_{x_2} \psi\right),
\end{equation}
where the domain is given by 
\begin{equation*}
	D(\mcA_k) = \left\{ u \in H_0^1(\La;w) \colon  \partial_{x_2}(x_2^\alpha \partial_{x_2} u) \in L^2(\La)  \right\}. 
\end{equation*}

We first analyze the baseline self-adjoint operator corresponding to the zero mode $k = 0$. Let $\mcA_0 = -\partial_{x_2}(x_2^\alpha \partial_{x_2})$. By Remark \ref{10.30.R2} (3), the embedding of $D(\mcA_0)$ into $L^2(\La)$ is compact. According to the classical spectral theorem for compact self-adjoint operators (or the Hilbert-Schmidt theorem), the spectrum $\sigma(\mcA_0)$ consists entirely of an isolated, discrete sequence of real positive eigenvalues:
\begin{equation}\label{10.31.13}
	0 < \lambda_1^{(0)} \le \lambda_2^{(0)} \le \dots \le \lambda_j^{(0)} \le \dots \to +\infty,
\end{equation}
and their corresponding eigenfunctions $\{\phi_j^{(0)}(x_2)\}_{j=1}^\infty$ form a complete orthonormal basis of $L^2(\La)$.

For later use, the zero-mode eigenvalues can be computed explicitly.
Set
\[
\kappa=\frac{2-\alpha}{2},\qquad
\nu=\frac{1-\alpha}{2-\alpha}.
\]
Substituting
$v(x_2)=x_2^{(1-\alpha)/2}
W(\sqrt{\lambda}\,x_2^\kappa/\kappa)$ into
$-(x_2^\alpha v')'=\lambda v$ gives Bessel's equation
\[
z^2W''+zW'+(z^2-\nu^2)W=0.
\]
The condition $v(0)=0$ selects the $J_\nu$ branch, and $v(1)=0$
therefore yields
\[
\lambda_n^{(0)}=\kappa^2j_{\nu,n}^2,
\qquad
j_{\nu,n}=\pi\left(n+\frac\nu2-\frac14\right)+O(n^{-1}).
\]
Here $j_{\nu,n}$ is the $n$th positive zero of $J_\nu$.
In particular, these eigenvalues are simple and
\[
\lambda_n^{(0)}\sim\kappa^2\pi^2n^2,
\qquad
\lambda_{n+1}^{(0)}-\lambda_n^{(0)}
\sim2\kappa^2\pi^2n.
\]

Now, we consider the operator $\mcA_k$ for the non-zero modes $k \neq 0$. We establish the following spectral properties in Lemma \ref{10.31.L3}. 
	
\begin{lemma}\label{10.31.L3}
		Let $\alpha \in (0,1)$, $\gamma \in \{1,2\}$, and $k \in \mathbb{Z}- \{0\}$. The localized non-self-adjoint operator $\mcA_k$ defined in \eqref{10.31.12} enjoys the following spectral properties:
		
		(1) The spectrum $\sigma(\mcA_k)$ consists entirely of isolated discrete complex eigenvalues with finite algebraic multiplicities, i.e., $\sigma(\mcA_k) = \{\la_n^{(k)}\}_{n \in \mathbb{N}^*}$, satisfying $\lim_{n \to \infty} |\la_n^{(k)}| = +\infty$.
		
		(2) The spectrum $\sigma(\mcA_k)$ is strictly contained in the right half of the complex plane and satisfies the uniform lower bound:
			\begin{equation}\label{10.31.15}
				\sigma(\mcA_k) \s  \left\{ z \in \mathbb{C} \colon \mathrm{Re}(z) \ge \frac{(1-\alpha)^2}{4} \right\}.
			\end{equation}  
\end{lemma} 

\begin{proof}
	We divide the proof into two main steps: establishing the compact resolvent property and estimating the numerical range.
	
	{\it Step 1}. Compactness of the resolvent.
	
	Throughout this proof, let $\mathcal H=L^2(\La;\mathbb C)$.
	Recall the zero-mode baseline operator $\mcA_0 = -\partial_{x_2}(x_2^\alpha \partial_{x_2})$, which is a self-adjoint, strictly positive closed operator on $L^2(\La)$ with a compact resolvent $\mcA_0^{-1}$ induced by the compact embedding $H_0^1(\La;w) \hookrightarrow L^2(\La)$ (see Remark \ref{10.30.R2}). 
	
	For any $k \neq 0$, the operator $\mcA_k$ can be viewed as a perturbation of $\mcA_0$:
	\begin{equation*}
		\mcA_k = \mcA_0 + \mathcal{B}_k,   \text{ where } \mathcal{B}_k \psi = -\im k x_2^\ga \psi.
	\end{equation*}
	Since the speed weight satisfies $x_2 \in \La$, the perturbation operator $\mathcal{B}_k$ is a bounded linear operator on $L^2(\La)$ with its operator norm bounded by $\|\mathcal{B}_k\|_{\mathcal{L}(L^2(\La))} \le |k|$. 
	
	Let $\lambda_1 > 0$ be the minimal eigenvalue of $\mcA_0$. Then $\la_1\geq \f{(1-\al)^2}{4}$ by Remark \ref{10.30.R2} (1).  Define the open connected domain $D = \mathbb{C}- [\lambda_1, +\infty)$, on which $(\mcA_0 - \lambda I)^{-1}$ is holomorphic and compact. We perform the algebraic factorization:
	\begin{equation}\label{11.01.1}
		\mcA_k - \lambda I = \left( I + \mathcal{B}_k(\mcA_0 - \lambda I)^{-1} \right) (\mcA_0 - \lambda I), \mbox{ for all } \lambda \in D.
	\end{equation}
	The operator-valued function $A(\lambda) \equiv -\mathcal{B}_k(\mcA_0 - \lambda I)^{-1}$ is compact on $L^2(\La)$ for each $\lambda \in D$ and depends holomorphically on $\lambda$. 
	
	Let $\lambda_0 = -r$ for a sufficiently large positive parameter $r > 0$. By applying the classic resolvent estimate for self-adjoint operators, the operator norm of the analytic compact family satisfies \begin{equation*}
		\|A(-r)\|_{\mathcal{L}(\mathcal{H})} \le \|\mcB_k\|\|(\mcA_0+r I)^{-1}\|\leq  \frac{|k|}{\lambda_1 + r}
	\end{equation*}
	by $\|(\mcA_0+rI)^{-1}\|=\sup_{\la\in \si(\mcA_0)}\f{1}{\la+r}$. 
	Choosing a fixed threshold $r \ge r_0 > |k| - \lambda_1$ rigorously ensures that $\|A(-r)\|_{\mathcal{L}(\mathcal{H})} < 1$. Consequently, by the Banach fixed point theorem, the operator $I - A(-r)$ is boundedly invertible via the absolutely convergent Neumann series, which locates an explicit basepoint of invertibility for $\mcA_k - \lambda_0 I$ by \eqref{11.01.1}.

	Let us define the ascending chain of closed subspaces as \begin{equation*}
		K_m = \ker\left(I - \frac{1}{\eta_n}R(\lambda_0; \mcA_k)\right)^m
	\end{equation*} 
	for each $m \in \N^*$, where $\eta_n$ is the eigenvalue of the compact operator  $R(\la_0;\mcA_k)=(\la_0I-\mcA_k)^{-1}$. It is clear from the algebraic structure that $K_1 \s  K_2 \s  \dots \s  K_m \s  K_{m+1} \s  \dots$. We claim that this chain must stabilize at a finite integer index $p \in \N^*$. 
	
	Suppose, ad absurdum, that the chain does not stabilize, which implies that the inclusion is strictly proper for all $m \ge 1$, i.e., $K_m \subsetneq K_{m+1}$. Since each $K_m$ is a closed subspace of the Hilbert space $\mathcal{H}$, by invoking the classical Riesz's geometric lemma, there exists a sequence of normalized generalized eigenfunctions $\{v_m\}_{m=1}^\infty$ such that $v_m \in K_{m+1}$, $\|v_m\|_{\mathcal{H}} = 1$, and it maintains a strict geometric distance from the preceding subspace:
	\begin{equation}\label{11.01.2}
		\mathrm{dist}(v_m, K_m) = \inf_{z \in K_m} \|v_m - z\|_{\mathcal{H}} \ge \frac{1}{2}, \mbox{ for all } m \ge 1.
	\end{equation}
	Now, let us analyze the behavior of this bounded sequence $\{v_m\}_{m=1}^\infty$ under the action of the compact resolvent operator $R \equiv R(\lambda_0; \mcA_k)$. For any integers $m > \ell \ge 1$, we evaluate the distance between their images under the identity perturbation operator $T \equiv I - \frac{1}{\eta_n}R$:
	\begin{equation}\label{11.01.3}
		\frac{1}{\eta_n} R v_m - \frac{1}{\eta_n} R v_\ell = v_m - (T v_m + v_\ell - T v_\ell).
	\end{equation}
	By utilizing the chain definition $K_m = \ker(T^m)$, we observe that:
	\begin{enumerate}
		\item Since $v_m \in K_{m+1} = \ker(T^{m+1})$, it follows that $T v_m \in \ker(T^m) = K_m$.
		\item Since $\ell < m$, we have $\ell + 1 \le m$, which implies $v_\ell \in K_{\ell+1} \subset K_m$.
		\item Similarly, $T v_\ell \in K_\ell \subset K_m$.
	\end{enumerate}
	Consequently, the entire grouped term inside the parenthesis in \eqref{11.01.3} belongs to the subspace $K_m$, i.e., $z_0 \equiv (T v_m + v_\ell - T v_\ell) \in K_m$. By substituting $z_0$ back into \eqref{11.01.3} and invoking the geometric distance constraint \eqref{11.01.2}, we yield:
	\begin{equation}\label{11.01.4}
		\left\| \frac{1}{\eta_n} R v_m - \frac{1}{\eta_n} R v_\ell \right\|_{\mathcal{H}} = \|v_m - z_0\|_{\mathcal{H}} \ge \mathrm{dist}(v_m, K_m) \ge \frac{1}{2}.
	\end{equation}
	Since $\eta_n \neq 0$ is fixed, \eqref{11.01.4} implies that $\|R v_m - R v_\ell\|_{\mathcal{H}} \ge \frac{|\eta_n|}{2} > 0$ for all $m > \ell \ge 1$. This indicates that the sequence $\{R v_m\}_{m=1}^\infty$ contains absolutely no Cauchy subsequences, meaning it cannot contain any strongly convergent subsequences in $\mathcal{H}$. 
	
	However, since $\{v_m\}_{m=1}^\infty$ is a bounded sequence ($\|v_m\|_{\mathcal{H}} = 1$) and $R$ is a compact operator on $\mathcal{H}$, the image sequence $\{R v_m\}_{m=1}^\infty$ must contain a strongly convergent subsequence. This directly yields a fundamental topological contradiction.

Thus the chain stabilizes: $K_p=K_{p+1}=\cdots$ for some finite $p$.
For each nonzero eigenvalue $\eta$ of the compact operator
$R=(\lambda_0I-\mcA_k)^{-1}$, the corresponding eigenvalue of
$\mcA_k$ is
\[
\lambda=\lambda_0-\eta^{-1}.
\]
Indeed, $Rv=\eta v$ implies
$(\lambda_0I-\mcA_k)v=\eta^{-1}v$.
The resolvent spectral correspondence identifies the associated
root space of $\mcA_k$ with $K_p$. By the Riesz--Schauder theorem,
this root space is finite-dimensional. Since $\mcA_k$ has compact
resolvent, all its finite spectral points are isolated eigenvalues
of finite algebraic multiplicity.

We also verify that there are infinitely many such eigenvalues.
Using the zero-mode eigenvalues computed above, set, for $n\ge2$,
\[
d_n=\frac13\min\{\lambda_n^{(0)}-\lambda_{n-1}^{(0)},
\lambda_{n+1}^{(0)}-\lambda_n^{(0)}\},
\qquad \Gamma_n=\{z:|z-\lambda_n^{(0)}|=d_n\}.
\]
The Bessel-zero asymptotics give $d_n\asymp n\to\infty$.
For fixed $k$ and sufficiently large $n$, the circles are disjoint
and, uniformly for $0\le\theta\le1$,
\[
\sup_{z\in\Gamma_n}
\|\theta\mcB_k(zI-\mcA_0)^{-1}\|
\le \frac{|k|}{d_n}<1.
\]
Consequently, $\Gamma_n\subset\rho(\mcA_0+\theta\mcB_k)$.
The spectral projections, with positively oriented contours,
\[
P_n(\theta)=\frac{1}{2\pi\mathrm{i}}
\int_{\Gamma_n}(zI-\mcA_0-\theta\mcB_k)^{-1}\,\mathrm{d}z
\]
depend continuously on $\theta$ in operator norm. Since $P_n(0)$
has rank one, each $P_n(\theta)$ has rank one.
Thus every sufficiently large circle encloses one eigenvalue of
$\mcA_k$, counted with algebraic multiplicity. In particular,
there are infinitely many eigenvalues. Enumerating the nonzero
eigenvalues $\eta_n$ of $R$ with algebraic multiplicities gives
\[
|\lambda_n^{(k)}|
=|\lambda_0-\eta_n^{-1}|\longrightarrow\infty,
\qquad \eta_n\longrightarrow0.
\]
This proves (1).

	{\it Step 2}. Numerical range and uniform spectral gap.
	
Let $\mu \in \sigma(\mcA_k)$ be an eigenvalue of $\mcA_k$, and let $\psi \in D(\mcA_k)$ be the corresponding normalized eigenfunction such that $\|\psi\|_{L^2(0,1)} = 1$ and $\mcA_k \psi = \mu \psi$. Taking the $L^2(\La)$ inner product of the eigenvalue equation with $\psi$, we obtain:
	\begin{equation}\label{11.01.5}
		\mu = \mu \|\psi\|_{L^2(\La)}^2 = \langle \mcA_k \psi, \psi \rangle_{L^2(\La)} = \int_0^1 \left( -\im k x_2^\ga \psi - \partial_{x_2}\left(x_2^\alpha \partial_{x_2} \psi\right) \right) \overline{\psi} \df x_2.
	\end{equation}
	Applying integration by parts to the second-order degenerate diffusion term, and utilizing the boundary conditions $\psi(0) = \psi(1) = 0$ and $\al\in (0,1)$, we have:
	\begin{equation*}
		-\int_0^1 [\partial_{x_2}(x_2^\alpha \partial_{x_2} \psi)] \overline{\psi} \df x_2 = \int_0^1 x_2^\alpha (\partial_{x_2} \psi)( \partial_{x_2} \overline{\psi}) \df x_2 = \int_0^1 x_2^\alpha |\partial_{x_2} \psi|^2 \df x_2.
	\end{equation*}
	Substituting this back into \eqref{11.01.5} yields:
	\begin{equation}\label{11.01.6}
		\mu = \int_0^1 x_2^\alpha |\partial_{x_2} \psi|^2 \df x_2 - \im k \int_0^1 x_2^\ga |\psi|^2 \df x_2.
	\end{equation}
Taking the real part of \eqref{11.01.6} eliminates the purely imaginary transport term and gives
	\begin{equation*}
		\Real \mu = \int_0^1 x_2^\alpha |\partial_{x_2} \psi|^2 \df x_2.
	\end{equation*}
	Finally, by applying the one-dimensional weighted Poincaré-Hardy inequality established in Lemma \ref{10.31.L1} to the normalized eigenfunction $\psi$, we deduce:
	\begin{equation*}
		\Real \mu = \int_0^1 x_2^\alpha |\partial_{x_2} \psi|^2 \df x_2 \geq \frac{(1-\alpha)^2}{4} \int_0^1 |\psi|^2 \df x_2 = \frac{(1-\alpha)^2}{4}.
	\end{equation*}
	Since $\mu \in \sigma(\mcA_k)$ was chosen arbitrarily, this proves that the entire spectrum is confined to the complex half-plane defined by $\Real z \geq \frac{(1-\alpha)^2}{4}$. The proof is now complete.
\end{proof}

\begin{lemma}[{\cite[Corollary 4.11, p.~923]{Shkalikov}}]\label{11.08.L1}
Let $H$ be a Hilbert space, and let $T$ be an operator on $H$ with compact resolvent and of finite order. Assume that $T$ is self-adjoint, or that $T$ is normal and its spectrum lies on a finite system of rays. If $B$ is an operator which is completely $T$ subordinated or is generated by a completely $T$ subordinated quadratic form, then the system of root vectors of $T+B$ is complete in $H$.
\end{lemma}

\begin{lemma}[{\cite[Theorem 6.9, p.~940]{Shkalikov}}]\label{11.08.L3}
Let $H$ be a Hilbert space, and let $T$ be a normal operator on $H$ of order $m>0$ with spectrum lying on a finite system of rays $\arg\lambda=\theta_j$, $j=1,\dots,N$. If an operator $B$ is generated by a form which is $p$-subordinate to $T$ for some $p\in [0,1)$ and $m^{-1}\le 1-p$, then the system of root vectors of $A=T+B$ forms an unconditional basis with parentheses in $H$. 
\end{lemma}

	We next establish the completeness needed for the density argument.
	\begin{lemma}\label{11.08.L2}
		Let $\alpha\in(0,1)$, $\gamma\in\{1,2\}$ and
		$k\in\mathbb Z\setminus\{0\}$.
		Choose a basis in each algebraic eigenspace of $\mcA_k$ and denote
		the resulting family of root vectors by $\{\phi_j^{(k)}\}_{j\ge1}$.
		Then
		\[
		\overline{\operatorname{span}\{\phi_j^{(k)}:j\ge1\}}^{L^2(\La)}
		=L^2(\La).
		\]
	\end{lemma}
	\begin{proof}
		The operator $\mcA_0$ is positive and self-adjoint.
		The Bessel formula above gives quadratic eigenvalue growth,
		so $\mcA_0$ has finite order $m=2$ and
		$\mcA_0^{-1}\in\mathcal S_q$ for every $q>1/2$.
		For fixed $k$, the perturbation $\mcB_kv=-\mathrm{i}kx_2^\gamma v$
		is bounded on $L^2(\La)$. In particular, for every $\varepsilon>0$,
		\[
		\|\mcB_kv\|_2\le |k|\|v\|_2
		\le\varepsilon\|\mcA_0v\|_2+|k|\|v\|_2,
		\qquad v\in D(\mcA_0).
		\]
		Thus $\mcB_k$ is completely $\mcA_0$-subordinated.
		Lemma~\ref{11.08.L1} proves the stated completeness.
	\end{proof}

		

		Let $E=H_0^1(\La;w)$ with the equivalent energy norm
		\[
		\|v\|_E^2=\int_0^1x_2^\alpha|v'(x_2)|^2\,\mathrm{d}x_2,
		\qquad \omega_0=\frac{(1-\alpha)^2}{4}.
		\]
		For each $k\in\mathbb Z$, let $S_k$ be the algebraic span of the
		root vectors of $\mcA_k$; for $k=0$ use its orthonormal eigenfamily.
		This span is dense in $L^2(\La)$ by Lemma~\ref{11.08.L2} and the
		self-adjoint spectral theorem for $k=0$.
		
		Fix $\lambda>0$ and write $R_k=(\lambda I+\mcA_k)^{-1}$.
		The spectral location proved above makes $R_k$ well-defined.
		Taking real parts in $(\lambda I+\mcA_k)v=f$ gives
		\[
		\lambda\|v\|_2^2+\|v\|_E^2
		=\operatorname{Re}(f,v),\qquad
		\|v\|_2\le\frac{\|f\|_2}{\lambda+\omega_0},\qquad
		\|R_kf\|_E^2\le\frac{\|f\|_2^2}{\lambda+\omega_0}.
		\]
		Hence $R_k:L^2(\La)\to E$ is continuous.
		Each finite-dimensional root space is invariant under $R_k$,
		and its restriction there is invertible. Therefore $R_kS_k=S_k$.
		Approximating $f\in L^2$ by elements of $S_k$ shows that $S_k$
		is dense in $R_kL^2=D(\mcA_k)$ in the energy topology.
		Finally, $C_c^\infty(0,1)\subset D(\mcA_k)$ is dense in $E$,
		so
		\[
		\overline{S_k}^{\|\cdot\|_E}=E.
		\]
		Only this density, rather than an energy-space Riesz basis,
		will be used below.

		By tensorizing these localized spectral components with the Fourier basis on the torus $\mathbb{T}_{x_1}$, we define the global generalized tensor product family $\mathcal{E}$ as follows:
		\begin{equation}\label{11.02.1}
			\mathcal{E} = \left\{ \Phi_j^{(k)}(x_1, x_2) \equiv \left(\frac{1}{\sqrt{2\pi}} e^{\mathrm{i} k x_1}\right) \phi_j^{(k)}(x_2) \colon k \in \mathbb{Z}, \, j \in \mathbb{N}^* \right\}.
		\end{equation}
		
		\begin{lemma}\label{10.31.L4}
			The linear span of the generalized tensor product family $\mathcal{E}$ defined in \eqref{11.02.1} is dense in $V_0$ with respect to the standard norm $\|\cdot\|_{V_0}$, i.e., 
			\begin{equation*}
				V_0 = \overline{\mathrm{span}\,\mathcal{E}}^{\|\cdot\|_{V_0}}. 
			\end{equation*}
		\end{lemma}

		\begin{proof}
			Let $\varphi\in V_0$ and $\varepsilon>0$, and put
			$e_k(x_1)=(2\pi)^{-1/2}e^{\mathrm{i}kx_1}$.
			By density of smooth functions in $V_0$ and horizontal Parseval,
			$\widehat\varphi_k\in E$ and
			\[
			\|\varphi\|_{V_0}^2
			=\sum_{k\in\mathbb Z}\|\widehat\varphi_k\|_E^2
			\]
			for the equivalent energy norms.
			Choose $K$ such that
			\[
			\sum_{|k|>K}\|\widehat\varphi_k\|_E^2<\varepsilon^2/2.
			\]
			For each of the finitely many indices $|k|\le K$, the energy
			density proved above gives $g_k\in S_k$ such that
			\[
			\sum_{|k|\le K}\|\widehat\varphi_k-g_k\|_E^2<\varepsilon^2/2.
			\]
			The function $g=\sum_{|k|\le K}e_kg_k$ belongs to
			$\operatorname{span}\mathcal E$. Again by Parseval,
			\[
			\|\varphi-g\|_{V_0}^2
			=\sum_{|k|>K}\|\widehat\varphi_k\|_E^2
			+\sum_{|k|\le K}\|\widehat\varphi_k-g_k\|_E^2
			<\varepsilon^2.
			\]
			This proves the required density without any uniform spectral
			projection bound in $k$.
		\end{proof}

		\begin{remark}\label{10.31.R3}
			Since $V_0$ is dense in $L^2(\Omega)$ and the embedding
			$V_0\hookrightarrow L^2(\Omega)$ is continuous,
			Lemma~\ref{10.31.L4} also gives
			\[
			\overline{\operatorname{span}\mathcal E}^{L^2(\Omega)}
			=L^2(\Omega).
			\]
			This is a completeness statement. No global Riesz basis property
			or spectral projection bound uniform in $k$ is required.
		\end{remark}

\subsection{Mild solution}

In this subsection, we present the mild solution to equation \eqref{10.31.1}.

\begin{theorem}\label{11.08.T1}
		Let $\alpha \in (0,1)$ and $\gamma \in \{1,2\}$. For any given initial data $\vp^0 \in L^2(\Om)$ and  $u \in  L^2(Q)$, the degenerate kinetic system \eqref{10.31.1} admits a unique mild solution $\vp \in C([0,T]; L^2(\Om))$. Moreover, this solution satisfies the uniform a priori estimate 
		\begin{equation}\label{11.03.1} 
			\|\vp(t)\|_{L^2(\Om)} \le e^{-\frac{(1-\alpha)^2}{4}t}\|\vp^0\|_{L^2(\Om)} + \int_0^t e^{-\frac{(1-\alpha)^2}{4}(t-s)}\|u(s)\|_{L^2(\Om)} \df s, \mbox{ for all } t \in [0,T]. 
		\end{equation} 
	\end{theorem}

\begin{proof}
	We work in the complexified spaces, with inner products linear in
	the first argument, and set
	\[
	q_\Omega(v)=\int_\Omega x_2^\alpha
	|\partial_{x_2}v|^2\,\mathrm{d}x,
	\qquad \omega_0=\frac{(1-\alpha)^2}{4}.
	\]
	The domain is dense because $C_c^\infty(\Omega)\subset D(\mcA)$.
	Put $e_k(x_1)=(2\pi)^{-1/2}e^{\mathrm{i}kx_1}$.
	For $v\in D(\mcA)$, let $v_N=\sum_{|k|\le N}e_k\widehat v_k$.
	Then $\widehat v_k\in D(\mcA_k)$ and
	$\widehat{\mcA v}_k=\mcA_k\widehat v_k$ in distributions.
	For these finite sums the separate integrations by parts give
	\begin{equation}\label{11.03.4}
		(\mcA v_N,v_N)_{L^2(\Omega)}
		=-\int_\Omega\left(x_2^\gamma\partial_{x_1}v_N
		+\partial_{x_2}(x_2^\alpha\partial_{x_2}v_N)\right)
		\overline{v_N}\,\mathrm{d}x.
	\end{equation}
	Periodicity gives
	\begin{equation}\label{11.03.2}
		\operatorname{Re}\int_\Omega x_2^\gamma
		\partial_{x_1}v_N\,\overline{v_N}\,\mathrm{d}x=0,
	\end{equation}
	and the vertical boundary conditions give
	\begin{equation}\label{11.03.3}
		-\operatorname{Re}\int_\Omega
		\partial_{x_2}(x_2^\alpha\partial_{x_2}v_N)
		\overline{v_N}\,\mathrm{d}x=q_\Omega(v_N).
	\end{equation}
	Since $v_N\to v$ in $V_0$ and
	$\mcA v_N\to\mcA v$ in $L^2$, passage to the limit yields
	\begin{equation}\label{11.03.5}
		\operatorname{Re}(\mcA v,v)_{L^2(\Omega)}
		=q_\Omega(v)\ge\omega_0\|v\|_{L^2(\Omega)}^2.
	\end{equation}
	
	To prove closedness, suppose that $v_n\to v$ and
	$\mcA v_n\to f$ in $L^2(\Omega)$.
	Applying the preceding identity to $v_n-v_m$ gives
	\[
	q_\Omega(v_n-v_m)
	\le\|\mcA(v_n-v_m)\|_2\|v_n-v_m\|_2\longrightarrow0.
	\]
	Thus $v_n$ converges in $V_0$ to $v$, so the boundary conditions
	are preserved. Distributional convergence then gives $\mcA v=f$.
	Hence $v\in D(\mcA)$ and $\mcA$ is closed.
	
	Next fix $\lambda>0$ and $f\in L^2(\Omega)$.
	For each $k$, the continuous sesquilinear form on $E$,
	\[
	a_{\lambda,k}(u,z)=\lambda(u,z)_{L^2(\La)}
	+\int_0^1x_2^\alpha u'\overline{z'}\,\mathrm{d}x_2
	-\mathrm{i}k\int_0^1x_2^\gamma u\overline z\,\mathrm{d}x_2,
	\]
	has real part $\lambda\|u\|_2^2+\|u\|_E^2$ on the diagonal.
	Lax--Milgram therefore produces a unique $u_k\in E$ with
	$a_{\lambda,k}(u_k,z)=(\widehat f_k,z)$ for all $z\in E$.
	Testing against compactly supported functions shows
	$\partial_{x_2}(x_2^\alpha u_k')\in L^2$; thus
	$u_k\in D(\mcA_k)$ and $(\lambda I+\mcA_k)u_k=\widehat f_k$.
	Taking real parts gives
	\[
	\|u_k\|_2\le\frac{\|\widehat f_k\|_2}{\lambda+\omega_0},
	\qquad
	\|u_k\|_E^2\le
	\frac{\|\widehat f_k\|_2^2}{\lambda+\omega_0}.
	\]
	Consequently $u=\sum_ke_ku_k$ converges in $V_0$, and
	\[
	\sum_k\|u_k\|_E^2\le
	\frac{\|f\|_{L^2(\Omega)}^2}{\lambda+\omega_0}<\infty.
	\]
	The modal equations imply $\mcA u=f-\lambda u\in L^2$ in
	distributions. Thus $u\in D(\mcA)$ and
	$(\lambda I+\mcA)u=f$.
	
	The Lumer--Phillips theorem, also applied to
	$-\mcA+\omega_0I$, now gives
	\begin{equation}\label{11.03.6}
		\|e^{-t\mcA}\|_{\mathcal L(L^2(\Omega))}
		\le e^{-\omega_0t},\qquad t\ge0.
	\end{equation}
	The unique mild solution is
	\[
	\vp(t)=e^{-t\mcA}\vp^0+
	\int_0^t e^{-(t-s)\mcA}u(s)\,\mathrm{d}s.
	\]
	It belongs to $C([0,T];L^2(\Omega))$, and the last estimate
	implies \eqref{11.03.1}.
\end{proof}

\subsection{Method of separation of variables: step 2}

From Lemma \ref{10.31.L4}, for the equation \eqref{10.31.1}, we denote the partial Fourier expansion with respect to the horizontal variable $x_1$ as
\begin{equation}\label{11.03.7}
	\vp(x_1, x_2, t)=\sum_{k\in\Z}\wh \vp_k(x_2, t)\left(\f{1}{\sqrt{2\pi}}e^{\im kx_1}\right).
\end{equation}
Then, by projecting the governing equation \eqref{10.31.1} onto each orthogonal Fourier mode, we obtain the following decoupled family of localized one-dimensional evolution problems:
\begin{equation}\label{10.03.8}
	\begin{cases}
		\pt_t\wh\vp_k - \im kx_2^\ga \wh\vp_k - \pt_{x_2}\left(x_2^\al \pt_{x_2}\wh\vp_k\right)=\wh u_k, &\mbox{in } \La\ts (0,T),\\[1mm]
		\wh\vp_k=0, &\mbox{on }\pt\La\ts (0,T),\\[1mm]
		\wh\vp_k(x_2, 0)=\wh\vp_k^0(x_2), &\mbox{in }\La
	\end{cases}
\end{equation}
for all $k\in\Z$, where the localized initial condition is defined via the $L^2(\mbT_{x_1})$ partial inner product:
\begin{equation*}
	\wh\vp_k^0(x_2) = \left(\vp^0(\cdot, x_2), \f{1}{\sqrt{2\pi}}e^{\im kx_1}\right)_{L^2(\mbT_{x_1})} = \f{1}{\sqrt{2\pi}}\int_{\mbT_{x_1}} \vp^0(x_1, x_2) e^{-\im k x_1} \df x_1.
\end{equation*}

It is obvious that equation \eqref{10.03.8} is a complex partial differential equation, and its solution $\wh\vp_k\in C([0,T];L^2(\La;\Cc))$. The existence and uniqueness of the localized mild solution are rigorously guaranteed by the localized semigroup theory. 

Indeed, for each fixed $k \in \mathbb{Z}$, the closed and densely defined operator $-\mcA_k$ forms a bounded perturbation of the self-adjoint baseline operator $-\mcA_0$. For every $\psi\in D(\mcA_k)$, integration by parts and
Remark~\ref{10.30.R2} give
\begin{equation}\label{11.03.9}
	\Real \langle \mcA_k \psi, \psi \rangle_{L^2(\La)} \ge \frac{(1-\alpha)^2}{4}\|\psi\|_{L^2(\La)}^2.
\end{equation}
Thus $-\mcA_k$ is strictly dissipative. We verify the range
condition below.
 
To rigorously verify the localized range condition $R(\lambda I + \mcA_k) = L^2(\Lambda)$ for $\lambda > 0$, we invoke the bounded perturbation framework developed in Lemma \ref{10.31.L3}. Let $\mcB_k = -\mathrm{i}k x_2^\gamma$ be the bounded coupling operator on $L^2(\Lambda)$ with $\|\mcB_k\|_{\mcL(L^2(\La))} \le |k|$. For any $\lambda > 0$, the strictly positive baseline self-adjoint operator $\mcA_0 = -\partial_{x_2}(x_2^\alpha \partial_{x_2})$ satisfies $\lambda \in \rho(-\mcA_0)$, implying that $(\lambda I + \mcA_0): D(\mcA_0) \to L^2(\Lambda)$ is a bijective mapping with a compact resolvent $(\lambda I + \mcA_0)^{-1}$. We perform the following algebraic splitting:
\begin{equation}\label{11.03.8}
	\lambda I + \mcA_k = \left[ I + \mcB_k (\lambda I + \mcA_0)^{-1} \right] (\lambda I + \mcA_0).
\end{equation}
The operator-valued function $T(\lambda) \equiv -\mcB_k (\lambda I + \mcA_0)^{-1}$ forms a holomorphic family of compact operators on $L^2(\Lambda)$ for $\text{Re}\,\lambda > 0$. According to the Analytic Fredholm Theorem, the bounded invertibility of the identity perturbation $I - T(\lambda)$ is governed by the Fredholm alternative. Since $\|T(r)\|_{\mathcal{L}(L^2(\La))} \le \f{|k|}{r+\la_1^{(0)}} < 1$ for a sufficiently large real parameter $r \gg 0$, the Neumann series guarantees the existence of a basepoint of invertibility. Furthermore, the strict numerical range estimate established via dissipativity in \eqref{11.03.9} (see also \eqref{10.31.15}) rigorously precludes the existence of any discrete eigenvalues of $-\mcA_k$  on the positive real axis $(0,+\infty)$, and by the Fredholm alternative, the operator $I - T(\lambda)$ is a surjective mapping onto $L^2(\Lambda)$. Combining this with the surjectivity of the self-adjoint baseline component $\lambda I + \mcA_0$, the algebraic factorization \eqref{11.03.8} systematically yields $R(\lambda I + \mcA_k) = L^2(\Lambda)$ for all $\lambda > 0$.

Consequently, by virtue of the Lumer-Phillips theorem \cite[Theorem 4.3, p.~14]{Pazy}, $-\mcA_k$ generates a strongly continuous semigroup $\{e^{-t\mcA_k}\}_{t \ge 0}$ on the complex Hilbert space $L^2(\La;\mathbb{C})$, which systematically ensures the well-posedness of the decoupled evolution system \eqref{10.03.8} via Duhamel's formula
\begin{equation*}
	\wh\vp_k(t) = e^{-t\mcA_k}\wh\vp_k^0 + \int_0^t e^{-(t-s)\mcA_k} \wh u_k(s) \df s. 
\end{equation*}

\begin{lemma}\label{11.05.L1}
	Let $k\in\Z$, $\wh\vp_k^0\in L^2(\La;\mathbb{C})$, and $\wh u_k=0$. Then, the unique mild solution of the homogeneous decoupled system \eqref{10.03.8} satisfies the following uniform exponential dissipation estimate:
	\begin{equation*}
		\|\wh\vp_k(t)\|_{L^2(\La)}\leq  \|\wh\vp_k^0\|_{L^2(\La)}e^{-2\omega_0 t}, \text{ for all } t \in [0,T],
	\end{equation*}
	where the positive decay rate $\omega_0 = \frac{(1-\alpha)^2}{4}$.
\end{lemma}

\begin{proof}
	For $v\in E$, the boundary condition $v(0)=0$ and
	Cauchy--Schwarz give
	\[
	|v(x_2)|^2\le\frac{x_2^{1-\alpha}}{1-\alpha}
	\int_0^1y^\alpha|v'(y)|^2\,\mathrm{d}y.
	\]
	Integrating with respect to $x_2$ yields
	\[
	\int_0^1x_2^\alpha|v'|^2\,\mathrm{d}x_2
	\ge(1-\alpha)(2-\alpha)\|v\|_2^2
	\ge2\omega_0\|v\|_2^2.
	\]
	First take $\widehat\vp_k^0\in D(\mcA_k)$.
	For the corresponding strong solution, the skew potential has
	zero real energy contribution. Therefore
	\begin{equation}\label{11.05.1}
		\frac12\frac{\mathrm d}{\mathrm dt}\|\widehat\vp_k(t)\|_2^2
		=\operatorname{Re}(-\mcA_k\widehat\vp_k(t),
		\widehat\vp_k(t))_{L^2}
		=-\int_0^1x_2^\alpha|\partial_{x_2}\widehat\vp_k|^2
		\,\mathrm{d}x_2.
	\end{equation}
	Thus
	\[
	\frac{\mathrm d}{\mathrm dt}\|\widehat\vp_k(t)\|_2^2
	\le-4\omega_0\|\widehat\vp_k(t)\|_2^2,
	\qquad
	\|\widehat\vp_k(t)\|_2
	\le e^{-2\omega_0t}\|\widehat\vp_k^0\|_2.
	\]
	Approximation of an arbitrary $L^2$ initial datum by elements of
	$D(\mcA_k)$ and continuity of the semigroup prove the same bound
	for all mild solutions. The bound is independent of $k$.
\end{proof}

\section{Carleman estimate}\label{S3}

In this section, all  functions are complex-valued functions. Denote
\begin{equation}\label{11.18.1}
	\La_T=\La\ts (0,T),\quad p=\f{1+\al}{2},\quad r=\f{3-\al}{2},\quad h=\f{3\al-1}{2}. 
\end{equation}
Then $0<\al<p$, $h<1$, and $h+p=2\al$, and $2p-\al=1$, and $r+p-2=0$. 

\begin{lemma}\label{11.18.L1}
	Let $u\in D(\mcA_0)$, where $\mcA_0$ is defined in \eqref{10.31.12}. Then 
	\begin{equation*}
		\int_\La x_2^h |\pt_{x_2}u|^2\df x_2\leq C\|u\|_{D(\mcA_0)}^2,\quad \int_\La x_2^{h-2}|u|^2\df x_2\leq \f{4}{(1-h)^2}\int_\La x_2^h|\pt_{x_2}u|^2\df x_2. 
	\end{equation*}
\end{lemma}

\begin{proof}
	Since $u\in D(\mcA_0)$, we have $u\in H_0^1(\La;w)$ and $\pt_{x_2}(x_2^\al \pt_{x_2}u)\in L^2(\La)$, which implies $x_2^\al\pt_{x_2}u\in H^1(\La)$. In particular, $x_2^\al \pt_{x_2}u$ is an absolutely continuous function on $\overline{\La}$ with $\|x_2^\al \pt_{x_2}u\|_{L^\iy(\La)}\leq C\|u\|_{D(\mcA_0)}$. Setting $p:=2\al-h$, the condition $p<1$ ensures that $x_2^{-p}$ is integrable on $\La=(0,1)$, and hence
	\begin{equation*}
		\begin{split}
			\int_\La x_2^h|\pt_{x_2}u|^2\df x_2 
			&=\int_\La \left|x_2^{ \al} \pt_{x_2}u\right|^2 x_2^{-p}\df x_2 \leq \|x_2^\al \pt_{x_2}u\|_{L^\iy(\La)}^2\int_\La x_2^{-p}\df x_2\\
			&=\f{1}{1-p} \|x_2^\al \pt_{x_2}u\|_{L^\iy(\La)}^2\leq C \|u\|_{D(\mcA_0)}^2,
		\end{split}
	\end{equation*}
	where the positive constant $C$ depends only on $\al$ and $h$.
	
Next, since $0<\alpha<1$, functions in $H_0^1(\La;w)$ have zero traces at both endpoints. As $u\in D(\mcA_0)\subset H_0^1(\La;w)$, we have $u(0)=u(1)=0$. Hence, for every $x_2\in\La$,
	\begin{equation}\label{11.19.1}
		|u(x_2)|\leq \int_0^{x_2}s^{-\al}|s^\al \pt_s u(s)|\df s\leq \f{x_2^{1-\al}}{1-\al}\|x_2^\al \pt_{x_2}u\|_{L^\iy(\La)} \leq Cx_2^{1-\al}\|u\|_{D(\mcA_0)}.
	\end{equation} 
	Since $h+1-2\al = 1-p > 0$, we obtain
	\begin{equation*}
		0\leq x_2^{h-1}|u(x_2)|^2 \leq C x_2^{h+1-2\al}\|u\|_{D(\mcA_0)}^2 \to 0  \text{ as } x_2\to 0^+,
	\end{equation*}
	and clearly $\lim_{x_2\to 1^-} x_2^{h-1}|u(x_2)|^2 = 0$. Furthermore, on $\La=(0,1)$ we have
	\begin{equation*}
		x_2^{h-2}|u(x_2)|^2\leq Cx_2^{h-2\al} \|u\|_{D(\mcA_0)}^2=Cx_2^{-p}\|u\|_{D(\mcA_0)}^2\in L^1(\La),
	\end{equation*}
	which guarantees that $\int_\La x_2^{h-2}|u|^2\df x_2 < +\iy$. Integrating by parts and applying the Cauchy-Schwarz inequality, we deduce that
	\begin{equation*}
		\begin{split}
			\int_\La x_2^{h-2}|u|^2\df x_2 
			&= \f{1}{h-1}\int_\La \pt_{x_2}(x_2^{h-1}) |u|^2\df x_2 
			= -\f{2}{h-1} \operatorname{Re}\int_\La x_2^{h-1} \ol{u}\,\pt_{x_2}u\df x_2 \\
			&\leq \f{2}{|1-h|} \left(\int_\La x_2^{h-2}|u|^2\df x_2\right)^{\f{1}{2}} \left(\int_\La x_2^h |\pt_{x_2}u|^2\df x_2\right)^{\f{1}{2}},
		\end{split}
	\end{equation*}
	which gives
	\begin{equation*}
		\int_\La x_2^{h-2}|u|^2\df x_2 \leq \f{4}{(1-h)^2}\int_\La x_2^h|\pt_{x_2}u|^2\df x_2.
	\end{equation*}
	We prove this lemma. 
\end{proof}

\subsection{Conjugate operator decomposition and regularity of the transformed terms}\label{S3.1}

Let $k\in\mathbb Z$, $\gamma\ge1$,
$\widehat\vp_k^0\in D(\mcA_k^2)$ and
$\widehat u_k\in H^2(0,T;D(\mcA_k^2))$.
The semigroup preserves $D(\mcA_k^2)$, and Duhamel's formula,
viewed in this graph space, gives
\[
\widehat\vp_k\in C([0,T];D(\mcA_k^2))
\cap C^1([0,T];D(\mcA_k)).
\]

Let 
\begin{equation}\label{11.04.9}
	\Theta=\f{1}{[t(T-t)]^4},\quad \eta=x_2^r, \quad \xi=\Theta (2-\eta). 
\end{equation}
Then $\xi(1,t)=\Theta(t)$ and $\Theta\leq \xi\leq 2\Theta$. 
Denote
\begin{equation}\label{11.04.10}
	\psi=e^{-s\xi}\wh\vp_k  \iff  \wh\vp_k=e^{s\xi}\psi. 
\end{equation}
The transformed function satisfies
\begin{equation}\label{11.03.10}
	\begin{gathered}
		\psi(0,t)=\psi(1,t)=0,\qquad 0<t<T,\\
		\lim_{t\to0^+,\,T^-}
		\left(\|\psi(\cdot,t)\|_{L^2(\La)}+
		\|x_2^{\alpha/2}\partial_{x_2}\psi(\cdot,t)\|_{L^2(\La)}
		\right)=0.
	\end{gathered}
\end{equation}

And 
\eqref{10.03.8} becomes 
\begin{equation*}
	\begin{split}
		e^{-s\xi}\wh u_k
		&=P_1\psi+P_2\psi+K\psi, 
	\end{split}
\end{equation*}
where 
\begin{equation*}
	\begin{split}
		P_1\psi
		&:=\sum_{i=1}^3P_{1i}\psi=-\pt_{x_2}\left(x_2^\al \pt_{x_2}\psi\right)+s\psi\xi_t-s^2\psi x_2^\al (\pt_{x_2}\xi)^2, \\
		P_2\psi
		&:=\sum_{i=1}^3P_{2i}\psi=\psi_t -2sx_2^\al (\pt_{x_2}\psi)(\pt_{x_2}\xi)-s\psi \pt_{x_2}\left(x_2^\al \pt_{x_2}\xi\right),\\
		K\psi
		&:=-\im kx_2^\ga \psi. 
	\end{split}
\end{equation*}

\begin{remark}\label{11.19.R1}
	We justify the integrations by parts for the strong data above.
	In this remark write $x=x_2$, $a=x^\alpha$,
	$g=e^{-s\xi}$, and $q_*=(1-\alpha)/2>0$.
	For $u\in D(\mcA_0)$ the flux $au'$ belongs to $H^1(0,1)$,
	so it is bounded, and $u(x)=O(x^{1-\alpha})$.
	Since $g'=sr\Theta x^{r-1}g$,
	\[
	\mcA_0(gu)=g\mcA_0u-2ag'u'-(ag')'u,
	\qquad ag'=O(x^p),\quad (ag')'=O(x^{p-1}).
	\]
	Both additional terms are $O(x^{p-\alpha})=O(x^{q_*})$.
	The product $gu$ preserves the zero boundary values and the
	weighted energy. Hence $gD(\mcA_0)\subset D(\mcA_0)$.
	For $\gamma\ge1$ the identity
	\[
	\mcA_0(x^\gamma u)=x^\gamma\mcA_0u
	-2\gamma x^{\gamma-1}(x^\alpha u')
	-\gamma(\alpha+\gamma-1)x^{\alpha+\gamma-2}u
	\]
	shows similarly that $x^\gamma D(\mcA_0)\subset D(\mcA_0)$.
	Since $D(\mcA_k)=D(\mcA_0)$, it follows that
	$D(\mcA_k^2)=D(\mcA_0^2)$.
	The same multiplier estimates apply to $g_t$.
	Consequently $\psi$ and $\psi_t$ belong to $D(\mcA_0)$ for each
	interior time, with graph-norm bounds uniform on compact time
	subintervals. We do not require $g$ to preserve $D(\mcA_0^2)$.
	
	Perform spatial integrations by parts first on $(\varepsilon,1)$,
	where the coefficients are smooth and the strong-data solution
	has the required local Sobolev regularity. Set
	$F=a\psi_x$ and $b=a\xi_x=-r\Theta x^p$.
	At a fixed interior time, $F\in H^1$,
	$\psi,\psi_t=O(x^{1-\alpha})$, and the potentially singular
	boundary expressions satisfy
	\[
	\begin{aligned}
		\xi_x|F|^2&=O(x^{q_*}),&
		Fb'\overline\psi&=O(x^{q_*}),\\
		ab''|\psi|^2&=O(x^{q_*}),&
		F\overline{\psi_t}&=O(x^{1-\alpha}).
	\end{aligned}
	\]
	They vanish as $\varepsilon\downarrow0$.
	The singular bulk terms are controlled by
	$x^h|\psi_x|^2$ and $x^{h-2}|\psi|^2$, both integrable by
	Lemma~\ref{11.18.L1}; their endpoint bounds are $O(x^{-p})$
	with $p<1$. All other spatial boundary terms contain vanishing
	factors $\psi$ or $\psi_t$. At $x=1$ the only surviving boundary
	term in the full cross-term identity is
	$2s\int_0^T\xi_x(1,t)|\psi_x(1,t)|^2\,\mathrm{d}t$.
	
	For time integrations, first work on $(\delta,T-\delta)$.
	The strong-data graph norms are bounded on $[0,T]$, while
	$e^{-s\Theta}$ dominates all powers of $\Theta$ and its time
	derivatives. Thus the temporal boundary terms vanish as
	$\delta\downarrow0$, and dominated convergence gives the full
	space-time identities. This proves \eqref{11.03.10} and justifies
	the subsequent calculations without inferring absolute
	continuity of a flux derivative from that of the flux itself.
\end{remark}

\subsection{Computation}\label{S3.2}

(1) From \eqref{11.03.10} we obtain 
\begin{equation*}
	\begin{split}
		&(P_1\psi, P_{21}\psi)_{L^2(\La_T)}+\ol{(P_1\psi, P_{21}\psi)_{L^2(\La_T)}}\\
		&=-s\iint_{\La_T}|\psi|^2\xi_{tt}\df x_2\df t+2s^2 \iint_{\La_T} x_2^\al (\pt_{x_2}\xi)(\pt_{x_2}\xi_t) |\psi|^2 \df x_2 \df t. 
	\end{split}
\end{equation*}

(2) From \eqref{11.03.10} we obtain 
\begin{equation*}
	\begin{split}
		&(P_{11}\psi, P_{22}\psi+P_{23}\psi)_{L^2(\La_T)}+\ol{(P_{11}\psi, P_{22}\psi+P_{23}\psi)_{L^2(\La_T)}}\\
		&=2s\int_0^T\left.\left(x_2^{2\al}|\pt_{x_2}\psi|^2\pt_{x_2}\xi\right)\right|_{x_2=0}^{x_2=1}\df t+s\iint_{\La_T} |\psi|^2\pt_{x_2}\left[x_2^\al \pt_{x_2x_2}\left(x_2^\al \pt_{x_2}\xi\right)\right]\df x_2\df t\\
		&\hspace{4.5mm}+2\al s\iint_{\La_T} x_2^{2\al-1}|\pt_{x_2}\psi|^2\pt_{x_2}\xi\df x_2\df t-4s\iint_{\La_T} x_2^\al |\pt_{x_2}\psi|^2\pt_{x_2}(x_2^\al \pt_{x_2}\xi)\df x_2\df t. 
	\end{split}
\end{equation*}

(3) From \eqref{11.03.10} we obtain 
\begin{equation*}
	\begin{split}
		&(P_{12}\psi, P_{22}\psi+P_{23}\psi)_{L^2(\La_T)}+\ol{(P_{12}\psi, P_{22}\psi+P_{23}\psi)_{L^2(\La_T)}}\\
		&=2s^2\iint_{\La_T} |\psi|^2(\pt_{x_2}\xi_t)(x_2^\al \pt_{x_2}\xi)\df x_2\df t. 
	\end{split}
\end{equation*}

(4) From \eqref{11.03.10} we obtain 
\begin{equation*}
	\begin{split}
		&(P_{13}\psi, P_{22}\psi+P_{23}\psi)_{L^2(\La_T)}+\ol{(P_{13}\psi, P_{22}\psi+P_{23}\psi)_{L^2(\La_T)}}\\
		&=-2s^3\iint_{\La_T} |\psi|^2\left(x_2^\al \pt_{x_2}\xi\right) \pt_{x_2}\left[x_2^\al (\pt_{x_2}\xi)^2\right]\df x_2\df t. 
	\end{split}
\end{equation*}

(5) From \eqref{11.03.10} we obtain 
\begin{equation*}
	\begin{split}
		&(P_1\psi, K\psi)_{L^2(\La_T)}+\ol{(P_1\psi, K\psi)_{L^2(\La_T)}}\\
		&=\im \ga k\iint_{\La_T} x_2^{\al+\ga-1}\ol \psi \pt_{x_2}\psi\df x_2\df t-\im \ga k\iint_{\La_T} x_2^{\al+\ga-1} \psi \pt_{x_2}\ol\psi\df x_2\df t. 
	\end{split}
\end{equation*}

\subsection{Estimate}\label{S3.3}

(6) We compute
\begin{equation*}
	\begin{split}
		|\Theta'|\leq 4T\Theta^\f{5}{4},\quad |\Theta\Theta'|\leq 4T\Theta^\f{9}{4},  \quad |\Theta''|\leq 20T^2\Theta^\f{3}{2}. 
	\end{split}
\end{equation*}

(7) We compute
\begin{equation*}
	\begin{split}
		\pt_{x_2}\eta=r x_2^{r-1}, \quad \pt_{x_2}\xi=-r\Theta x_2^{r-1}, \quad \pt_{x_2x_2}\xi=-r(r-1)\Theta x_2^{r-2}, 
	\end{split}
\end{equation*}
and 
\begin{equation*}
	\begin{split}
		x_2^\al\pt_{x_2}\xi=-r\Theta x_2^p,\quad \pt_{x_2}\left(x_2^\al \pt_{x_2}\xi\right)=-rp\Theta x_2^{p-1}, \quad \pt_{x_2}\xi_t=-r\Theta'x_2^{r-1}, 
	\end{split}
\end{equation*}
and 
\begin{equation*}
	\begin{split}
		\pt_{x_2}\left[x_2^\al (\pt_{x_2}\xi)^2\right]&=r^2\Theta^2, \quad  \xi_{tt}=\Theta''(2-x_2^r)\\
		\pt_{x_2}\left[x_2^\al \pt_{x_2x_2}(x_2^\al \pt_{x_2}\xi)\right]
		&=-rp(p-1)(h-1)\Theta x_2^{h-2}.
	\end{split}
\end{equation*}

(8) Since $p-\al = \f{1-\al}{2} > 0$ and $x_2^\al \pt_{x_2}\psi \in H^1(\La) \s  L^\iy(\La)$, we have
\begin{equation*}
	x_2^{2\al}(\pt_{x_2}\xi)|\pt_{x_2}\psi|^2=-r\Theta x_2^{p-\al}|x_2^\al\pt_{x_2}\psi|^2=0  \mbox{ at }x_2=0.
\end{equation*}
Together with this and (1)-(5) and (7), we obtain 
\begin{equation*}
	\begin{split}
		&2\Real (P_1\psi, P_2\psi+K\psi)_{L^2(\La_T)}\\
		&=-2rs\int_0^T\Theta |\pt_{x_2}\psi(1,t)|^2\df t\\
		&\hspace{4.5mm}+2rs\iint_{\La_T} \Theta x_2^h|\pt_{x_2}\psi|^2\df x_2\df t+2r^3s^3\iint_{\La_T} \Theta^3x_2^p|\psi|^2\df x_2\df t\\
		&\hspace{4.5mm}-s\iint_{\La_T}\Theta'' (2-x_2^r)|\psi|^2\df x_2\df t+4r^2s^2\iint_{\La_T} \Theta\Theta' x_2|\psi|^2\df x_2\df t\\
		&\hspace{4.5mm}-rp(p-1)(h-1)s\iint_{\La_T} \Theta x_2^{h-2}|\psi|^2\df x_2\df t-2\ga k\Image \iint_{\La_T} x_2^{\al+\ga-1}\ol\psi \pt_{x_2}\psi\df x_2\df t\\
		&=:\sum_{i=1}^7R_i.
	\end{split}
\end{equation*}
Note that $rp(p-1)(h-1)=\f{3(1-\al)^2(1+\al)(3-\al)}{16}>0$. 

(9) From Lemma \ref{11.18.L1} and 
\begin{equation*}
	2r - rp(p-1)(h-1)\f{4}{(1-h)^2} = \f{2(2-\al)}{3}r = \f{(2-\al)(3-\al)}{3} > 0,
\end{equation*}
we obtain 
\begin{equation*}
	\begin{split}
		R_2+R_6
		&\geq \f{(2-\al)(3-\al)}{3}s\iint_{\La_T} \Theta x_2^h|\pt_{x_2}\psi|^2\df x_2\df t.
	\end{split}
\end{equation*}

From Lemma \ref{11.18.L1}, for all $u\in D(\mcA_0)$, 
\begin{equation}\label{11.19.2}
	\begin{split}
		\int_\La |u|^2\df x_2
		&\leq \int_\La x_2^{\al-1}|u|^2\df x_2\leq \left(\int_\La x_2^p|u|^2\df x_2\right)^\f{1}{2}\left(\int_\La x_2^{h-2}|u|^2\df x_2\right)^\f{1}{2}\\
		&\leq \f{2}{1-h}\left(\int_\La x_2^p|u|^2\df x_2\right)^\f{1}{2}\left(\int_\La x_2^h|\pt_{x_2}u|^2\df x_2\right)^\f{1}{2}.
	\end{split}
\end{equation}
Together with this and (6), we obtain 
\begin{equation*}
	\begin{split} 
		|R_4|
		&\leq 40T^2s\iint_{\La_T} \Theta^\f{3}{2} |\psi|^2\df x_2\df t\\
		&\leq \f{5T^6}{1-h}s^{-1}\int_0^T\left[s^\f{3}{2}\Theta^\f{3}{2}\left(\int_\La x_2^p|\psi|^2\df x_2\right)^\f{1}{2}\right]\left[s^\f{1}{2}\Theta^\f{1}{2}\left(\int_\La x_2^h|\pt_{x_2}\psi|^2\df x_2\right)^\f{1}{2}\right]\df t\\
		&\leq \f{5T^6}{2(1-h)}s^2\iint_{\La_T}\Theta^3x_2^p|\psi|^2\df x_2\df t+\f{5T^6}{2(1-h)}\iint_{\La_T}\Theta x_2^h|\pt_{x_2}\psi|^2\df x_2\df t; 
	\end{split}
\end{equation*}
Furthermore, since $p<1$ and $\Theta^{-\f{3}{4}}\leq \f{T^6}{64}$, we have
\begin{equation*}
	\begin{split} 
		|R_5|
		&\leq 16Tr^2s^2\iint_{\La_T} \Theta^\f{9}{4}x_2|\psi|^2\df x_2\df t\leq \f{1}{4}T^7r^2s^2\iint_{\La_T}\Theta^3 x_2^p|\psi|^2\df x_2\df t; 
	\end{split}
\end{equation*}
Furthermore, since $\ga\geq 1$ and $\f{p+h}{2}=\al$ (by \eqref{11.18.1}), we have $x_2^{\al+\ga-1}\leq x_2^{\f{p+h}{2}}$ on $\La$. Combining this with $\Theta^{-2}\leq \f{T^{16}}{256^2}$ and Young's inequality yields
\begin{equation*}
	\begin{split} 
		|R_7|
		&\leq \f{2}{256^2}|\ga| |k|T^{16}s^{-2}\iint_{\La_T} \left(s^\f{3}{2} \Theta^\f{3}{2}x_2^\f{p}{2} |\psi|\right)\left(s^\f{1}{2} \Theta^\f{1}{2}x_2^\f{h}{2}|\pt_{x_2}\psi|\right)\df x_2\df t\\
		&\leq \f{1}{256^2}|\ga||k|T^{16}s\iint_{\La_T} \Theta^3 x_2^p|\psi|^2\df x_2\df t+\f{1}{256^2}|\ga||k|T^{16}s^{-1}\iint_{\La_T} \Theta x_2^h|\pt_{x_2}\psi|^2\df x_2\df t. 
	\end{split}
\end{equation*}

(10) From (9) and the estimates for $R_4$, $R_5$, and $R_7$, there exists a constant $\la_0\geq 1$, depending only on $\al$ and $\ga$, such that for all $\la\geq \la_0$ and $s=(1+|k|^\f{1}{2})(1+T^6+T^7+T^8)\la$, 
\begin{equation*}
	\begin{split}
		&s\iint_{\La_T}\Theta x_2^h|\pt_{x_2}\psi|^2\df x_2\df t+s^3\iint_{\La_T}\Theta^3 x_2^p|\psi|^2\df x_2\df t\\
		&\leq C\|e^{-s\xi}\wh u_k\|_{L^2(\La_T)}^2+Cs\int_0^T\Theta |\pt_{x_2}\psi(1,t)|^2\df t, 
	\end{split}
\end{equation*}
where the positive constant $C$ depends only on $\al$.

Note that from \eqref{11.19.2}, we deduce that 
\begin{equation*}
	\begin{split}
		&s^2\iint_{\La_T} \Theta^2 |\psi|^2\df x_2\df t\leq \f{1}{1-h}s^3\iint_{\La_T}\Theta^3x_2^p|\psi|^2\df x_2\df t+\f{1}{1-h}s\iint_{\La_T} \Theta x_2^h |\pt_{x_2}\psi|^2\df x_2\df t. 
	\end{split}
\end{equation*}
Hence, there exists a constant $\la_0\geq 1$, depending only on $\al$ and $\ga$, such that for all $\la\geq \la_0$ and $s=(1+|k|^\f{1}{2})(1+T^6+T^7+T^8)\la$, 
\begin{equation}\label{11.19.3}
	\begin{split}
		&s\iint_{\La_T}\Theta x_2^h|\pt_{x_2}\psi|^2\df x_2\df t+s^2\iint_{\La_T} \Theta^2 |\psi|^2\df x_2\df t+s^3\iint_{\La_T}\Theta^3 x_2^p|\psi|^2\df x_2\df t\\
		&\leq C\|e^{-s\xi}\wh u_k\|_{L^2(\La_T)}^2+Cs\int_0^T\Theta |\pt_{x_2}\psi(1,t)|^2\df t, 
	\end{split}
\end{equation}
where the positive constants $C$ depend only on $\al$ and $\ga$. 

(11) Since 
\begin{equation*}
	\pt_{x_2}\wh\vp_k=\pt_{x_2}\left(e^{s\xi}\psi\right)=e^{s\xi}\left(\pt_{x_2}\psi+s\psi\pt_{x_2}\xi\right),
\end{equation*}
and $x_2^h(\pt_{x_2}\xi)^2 = r^2\Theta^2 x_2^{h+2r-2} = r^2\Theta^2 x_2^p$, we obtain
\begin{equation*}
	\begin{split}
		s\iint_{\La_T} \Theta x_2^h|\pt_{x_2}\wh\vp_k|^2e^{-2s\xi}\df x_2\df t
		&\leq 2s\iint_{\La_T}\Theta x_2^h|\pt_{x_2}\psi|^2\df x_2\df t + 2s^3\iint_{\La_T}\Theta x_2^h(\pt_{x_2}\xi)^2|\psi|^2\df x_2\df t\\
		&= 2s\iint_{\La_T}\Theta x_2^h|\pt_{x_2}\psi|^2\df x_2\df t + 2r^2 s^3\iint_{\La_T}\Theta^3 x_2^p|\psi|^2\df x_2\df t.
	\end{split}
\end{equation*}
Moreover, since $\psi(1,t)=0$ and $\xi(1,t)=\Theta(t)$, we have $\pt_{x_2}\psi(1,t) = e^{-s\Theta}\pt_{x_2}\wh\vp_k(1,t)$. Combining these relations with \eqref{11.19.3} and $\psi = e^{-s\xi}\wh\vp_k$, we deduce that
\begin{equation*}
	\begin{split}
		&s\iint_{\La_T}\Theta x_2^h|\pt_{x_2}\wh\vp_k|^2e^{-2s\xi}\df x_2\df t + s^2\iint_{\La_T}\Theta^2|\wh\vp_k|^2e^{-2s\xi}\df x_2\df t + s^3\iint_{\La_T}\Theta^3 x_2^p|\wh\vp_k|^2e^{-2s\xi}\df x_2\df t\\
		&\leq C\|e^{-s\xi}\wh u_k\|_{L^2(\La_T)}^2 + Cs\int_0^T \Theta |\pt_{x_2}\wh\vp_k(1,t)|^2e^{-2s\Theta}\df t,
	\end{split}
\end{equation*}
where the positive constants $C$ depend only on $\al$ and $\ga$.  

Since $1+T^6+T^7+T^8\le3(1+T^8)$ for all $T>0$, this proves the following proposition.

\begin{proposition}\label{11.04.P1}  
	Let $\alpha \in (0,1)$, $\gamma \in \mathbb{N}^*$, $\wh\vp_k^0\in D(\mcA_k^2)$ and $\wh u_k\in H^2(0,T; D(\mcA_k^2))$, and $k \in \mathbb{Z}$. Let $\Theta(t)$, $\xi(x_2, t)$, $p$, $r$, and $h$ be defined as in \eqref{11.18.1}-\eqref{11.04.9}. 
Then there exists $\lambda_0=\lambda_0(\alpha,\gamma)\ge1$
	such that, for all $\lambda\ge\lambda_0$ and the corresponding
	Carleman parameter

	\begin{equation}\label{11.04.6}
		s = (1+|k|^\frac{1}{2})(1+T^8) \lambda,
	\end{equation}
	the following parameter-weighted Carleman estimate holds uniformly for every solution $\wh\vp_k \in C^1([0,T]; D(\mcA_k))$ of the decoupled one-dimensional degenerate system \eqref{10.03.8}:
	\begin{equation}\label{11.04.11}
		\begin{split}
			&s\iint_{\La_T} \Theta x_2^h |\pt_{x_2}\wh\vp_k|^2 e^{-2s\xi}\df x_2\df t \\
			&\quad + s^2\iint_{\La_T}\Theta^2 |\wh\vp_k|^2 e^{-2s\xi}\df x_2\df t + s^3\iint_{\La_T} \Theta^3 x_2^p |\wh\vp_k|^2 e^{-2s\xi}\df x_2\df t\\
			&\leq Cs\int_0^T\Theta |\pt_{x_2}\wh\vp_k(1,t)|^2 e^{-2s\Theta}\df t + C\|e^{-s\xi}\wh u_k\|_{L^2(\La_T)}^2,
		\end{split}
	\end{equation}
where $C=C(\alpha,\gamma)>0$ is independent of $T$, $k$ and
	$\lambda$. For $\gamma=1$, the constants $C$ and $\lambda_0$
	depend only on $\alpha$.
\end{proposition}

\begin{remark}\label{11.04.R1}
	Although the derivation borrows core concepts from classical Carleman estimates for degenerate parabolic equations, a crucial distinction lies in the coupling of the high-frequency Fourier parameter $k$ arising from the skew-symmetric transport operator $K$ (see Step \textup{(5)}). Specifically, the modal scaling rule \eqref{11.04.6} is custom-tailored to absorb the non-self-adjoint transport interaction terms via the heightened powers of the time weight $\Theta$.
\end{remark}

\section{Observability}\label{S4}

\subsection{A uniform a priori estimate}

\begin{lemma}\label{11.05.L2}
	Let $\ga\in \{1,2\}$, $k \in \Z- \{0\}$, $\wh\vp_k^0 \in L^2(\La;\Cc)$, and $\wh u_k=0$. Then, the unique mild solution $\wh\vp_k$ of the homogeneous decoupled system \eqref{10.03.8} satisfies
	\begin{equation}\label{11.21.2}
		\|\wh\vp_k(t)\|_{L^2(\La)} \leq C \|\wh\vp_k^0\|_{L^2(\La)} e^{-c_0 |k|^\be t}, \text{ for all } t \in [0,T],
	\end{equation}
	where $C \geq 1$ is a positive constant, $c_0 > 0$ is a dissipation rate constant, and the spectral enhancement exponent $\be$ is explicitly given by
	\begin{equation}\label{11.05.2}
		\be = \frac{2-\al}{2+\ga-\al} \in (0,1).
	\end{equation}
	Both positive constants $C$ and $c_0$ depend exclusively on $\al$ and $\ga$, remaining uniform with respect to $k \in \Z- \{0\}$ and $t\in [0,T]$.
\end{lemma}

\begin{proof}
	We divide the proof into five steps. 
	
	{\it Step 1}. Rescaling transformation.
	
	Define
	\begin{equation*}
		\de=|k|^{-\f{1}{2+\ga-\al}},\quad R=\de^{-1}, 
	\end{equation*}
	and 
	\begin{equation*}
		\begin{split}
			U_\de: L^2(\La)\ra L^2(0,R), \quad [U_\de(u)](y_2)=\de^\f{1}{2}u(\de y_2). 
		\end{split}
	\end{equation*}
	Then $U_\de$ is a unitary operator with inverse $[U_\de^{-1}(v)](x_2)=\de^{-\f{1}{2}}v(\de^{-1}x_2)$ for all $v\in L^2(0,R)$, and 
	\begin{equation*}
		\begin{split}
			(U_\de \mcA_k U_{\de}^{-1})(v)=-\de^{\al-2}\pt_{y_2}\left(y_2^\al \pt_{y_2}v(y_2)\right)-\im k \de^\ga y_2^\ga v(y_2). 
		\end{split}
	\end{equation*}
	Denote 
	\begin{equation}\label{11.20.3}
		\mcB_R v=-\pt_{y_2}(y_2^\al \pt_{y_2}v)-\im (\sgn k)y_2^\ga v. 
	\end{equation}
	Since $\de^{\al-2} = |k|\de^\ga = |k|^\be$, we obtain  
	\begin{equation}\label{11.20.1}
		\begin{split}
			(U_\de\mcA_kU_\de^{-1})(v)=|k|^\be \mcB_Rv, 
		\end{split}
	\end{equation}
	where $\be$ is defined in \eqref{11.05.2}. 
	
	{\it Step 2}. Maximal monotonicity and semigroup of $\mcB_R$.
	
	We define the domain of $\mcB_R$ by 
	\begin{equation}\label{11.20.2}
		D(\mcB_R)=\left\{u\in H_0^1((0,R); y_2^\al)\colon \mcB_R u\in L^2(0,R)\right\}, 
	\end{equation}
	where the inner product and norm on $H_0^1((0,R);y_2^\al)$ are defined (analogously to Section \ref{S2.2}) by
	\begin{equation*}
		(u,v)_{H_0^1((0,R);y_2^\al)}=\int_{(0,R)} y_2^\al (\pt_{y_2}u)(\pt_{y_2}\ol v)\df y_2,\quad \|u\|_{H_0^1((0,R);y_2^\al)}=(u,u)_{H_0^1((0,R); y_2^\al)}^\f{1}{2}. 
	\end{equation*}
	For each $u, v\in H_0^1((0,R);y_2^\al)$, we define the sesquilinear form
	\begin{equation}\label{11.20.4}
		b_R: H_0^1((0,R);y_2^\al)\ts H_0^1((0,R); y_2^\al)\ra \Cc
	\end{equation}
	by  
	\begin{equation*}
		b_R(u,v)=\int_{(0,R)} y_2^\al (\pt_{y_2}u)(\pt_{y_2}\ol v)\df y_2-\im (\sgn k)\int_{(0,R)}y_2^\ga u\ol v\df y_2. 
	\end{equation*}
	For each fixed $R>0$, by the weighted Hardy-Poincar\'e inequality on $(0,R)$, $b_R$ is bounded (continuous) on $H_0^1((0,R);y_2^\al)$. Moreover, for each $\la>0$, we have  
	\begin{equation*}
		\begin{split}
			\Real\left(b_R(u,u)+\la \|u\|_{L^2(0,R)}^2\right)=\int_{(0,R)}y_2^\al |\pt_{y_2}u|^2\df y_2+\la \int_{(0,R)}|u|^2\df y_2.
		\end{split}
	\end{equation*} 
	By the Lax-Milgram theorem, for each $f\in L^2(0,R)$, there exists a unique $u\in H_0^1((0,R);y_2^\al)$ such that
	\begin{equation*}
		b_R(u,v) + \la (u,v)_{L^2(0,R)} = (f,v)_{L^2(0,R)},  \text{ for all } v\in H_0^1((0,R);y_2^\al).
	\end{equation*}
	Testing with $v=\phi\in C_0^\iy(0,R)$ yields $-\pt_{y_2}(y_2^\al\pt_{y_2}u) - \im(\sgn k)y_2^\ga u = f - \la u$ in $\mathcal{D}'(0,R)$. Since $f - \la u \in L^2(0,R)$, we deduce that $\mcB_R u \in L^2(0,R)$, and thus $u\in D(\mcB_R)$ with $(\mcB_R+\la I)u = f$. Consequently, $\Ran(\mcB_R+\la I)=L^2(0,R)$. 
	Since $C_0^\iy(0,R)\s D(\mcB_R)$ is dense in $L^2(0,R)$, $\mcB_R$ is a maximal monotone (m-accretive) operator on $L^2(0,R)$. Hence, $-\mcB_R$ generates a contraction semigroup $S_R(t)=e^{-t\mcB_R}$ on $L^2(0,R)$ (see \cite[Theorem 4.3 in Chapter 1, p.~14]{Pazy}).
	
	{\it Step 3}. Energy and pointwise estimates. 
	
	We shall show 
\begin{equation}\label{11.20.5}
\begin{split}
	\left\|(\mcB_R+\im\mu)v\right\|_{L^2(0,R)}\geq \ka \|v\|_{L^2(0,R)}, \mbox{ for all } \mu\in\R \mbox{ and } v\in D(\mcB_R), 
\end{split}
\end{equation}
where the positive constant $\ka$ depends only on $\al$. It suffices to consider $\sgn k=1$. The case $\sgn k=-1$ follows by complex conjugation and the	replacement of $\mu$ by $-\mu$.
	
	Let $v\in D(\mcB_R)$ with $\|v\|_{L^2(0,R)}=1$. Write 
	\begin{equation}\label{11.20.7}
		f=(\mcB_R+\im \mu)v=-\pt_{y_2}(y_2^\al \pt_{y_2} v)-\im (y_2^\ga -\mu)v.
	\end{equation}
Multiplying both sides by $\overline v$, integrating over $(0,R)$, and taking real parts, we obtain
	\begin{equation}\label{11.20.6}
		\int_{(0,R)} y_2^\al |\pt_{y_2}v|^2\df y_2=\Real (f,v)_{L^2(0,R)}\leq \|f\|_{L^2(0,R)}. 
	\end{equation}
	
Let $\chi$ be a real  bounded  Lipschitz function on $(0,R)$. Then $\chi v\in H_0^1((0,R);y_2^\al)$. Testing \eqref{11.20.7} with $\chi v$ and taking the imaginary part, we obtain
	\begin{equation}\label{11.20.8}
		\begin{split}
			\int_{(0,R)}\chi(y_2^\ga -\mu) |v|^2\df y_2=-\Image \int_{(0,R)} \chi f\ol v\df y_2+\Image \int_{(0,R)} y_2^\al (\pt_{y_2}\chi)(\pt_{y_2}v)\ol v\df y_2. 
		\end{split}
	\end{equation}
	
	Now, for $\mu\leq 2^\ga$, taking 
	\begin{equation*}
		\chi(y_2)=
		\begin{cases}
			0, &y_2\leq 2,\\
			y_2-2, &2<y_2<3,\\
			1, &y_2\geq 3.
		\end{cases}
	\end{equation*}
On the support of $\chi$, the factor $y_2^\ga -\mu$ is nonnegative, and $y_2^\ga-\mu$ is at least $1$ for $y_2\geq 3$. From \eqref{11.20.8}, we get
\begin{equation}\label{11.20.9}
\begin{split}
\int_{(3,R)}|v|^2\df y_2
&\leq \|f\|_{L^2(0,R)}+\|v\|_{H_0^1((0,R);y_2^\al)}\left(                 \int_{(2,3)\cap(0,R)}y_2^\al |v|^2\df y_2       \right)^\f{1}{2}\\
&\leq \|f\|_{L^2(0,R)}+3^\f{\al}{2}\|v\|_{H_0^1((0,R);y_2^\al)}. 
\end{split}
\end{equation}
	Note that from
	\begin{equation*}
		\begin{split}
			|v(y_2)|^2\leq \left|\int_0^{y_2} (\pt_{y_2} v)(s)\df s\right|^2\leq \f{y_2^{1-\al}}{1-\al}\int_{(0,R)}y_2^\al |\pt_{y_2}v|^2\df y_2
		\end{split}
	\end{equation*} 
	we obtain 
	\begin{equation*}
		\begin{split}
			\int_0^{\min\{3,R\}}|v|^2\df y_2\leq \f{3^{2-\al}}{(1-\al)(2-\al)}\int_{(0,R)}y_2^\al |\pt_{y_2}v|^2\df y_2. 
		\end{split}
	\end{equation*}
	Combining \eqref{11.20.6} and \eqref{11.20.9}, we get 
	\begin{equation}\label{11.20.11}
		\begin{split}
			1
			&\leq \|f\|_{L^2(0,R)}+3^\f{\al}{2}\|v\|_{H_0^1((0,R); y_2^\al)}+\f{3^{2-\al}}{(1-\al)(2-\al)}\|v\|_{H_0^1((0,R); y_2^\al)}^2\\
			&\leq 3^\f{\al}{2}\|f\|_{L^2(0,R)}^\f{1}{2}+\left(1+\f{3^{2-\al}}{(1-\al)(2-\al)}\right)\|f\|_{L^2(0,R)}. 
		\end{split}
	\end{equation}
	
	For $\mu>2^\ga$, taking 
	\begin{equation*}
		\chi(y_2)=
		\begin{cases}
			-1, &y_2\leq \mu^\f{1}{\ga}-\f{1}{4}\mu^\f{\al}{2\ga}, \\
			4\mu^{- \f{\al}{2\ga}} (y_2-\mu^\f{1}{\ga}), &|y_2-\mu^\f{1}{\ga}|<\f{1}{4}\mu^\f{\al}{2\ga}, \\
			1, &y_2\geq \mu^\f{1}{\ga}+\f{1}{4}\mu^\f{\al}{2\ga}, 
		\end{cases}
	\end{equation*}
	then $\chi(y_2)(y_2^\ga-\mu)\geq 0$. For each $y_2\notin I:=(\mu^\f{1}{\ga}-\f{1}{4}\mu^\f{\al}{2\ga}, \mu^\f{1}{\ga}+\f{1}{4}\mu^\f{\al}{2\ga})\cap (0,R)$ and $\ga\in \{1,2\}$, we have 
	\begin{equation*}
		\begin{split}
			|y_2^\ga -\mu|\geq |y_2-\mu^\f{1}{\ga}|\geq \f{1}{4}\mu^\f{\al}{2\ga}\geq \f{1}{4}. 
		\end{split}
	\end{equation*}
	For each $y_2\in I$, we have (note that $\f{1}{4}\mu^\f{\al}{2\ga}<\f{1}{2}\mu^\f{1}{\ga}$)
	\begin{equation*}
		\f{1}{2} \mu^\f{1}{\ga}< y_2< \f{3}{2}\mu^\f{1}{\ga}, \mbox{ and } y_2^\al |\pt_{y_2}\chi|^2\leq 16\left(\f{3}{2}\right)^\al. 
	\end{equation*}
	Therefore, from \eqref{11.20.8}, we deduce that 
	\begin{equation}\label{11.20.10}
		\begin{split}
			\int_{(0,R)-I}|v|^2\df y_2
			&\leq 4\|f\|_{L^2(0,R)}+4\|v\|_{H_0^1((0,R);y_2^\al)}\left(\int_{(0,R)}y_2^\al |\pt_{y_2}\chi|^2 |v|^2\df y_2\right)^\f{1}{2}\\
			&\leq 4\|f\|_{L^2(0,R)}+16\left(\f{3}{2}\right)^\f{\al}{2} \|v\|_{H_0^1((0,R); y_2^\al)}. 
		\end{split}
	\end{equation}
	And from 
	\begin{equation*}
		\begin{split}
			|v(y_2)|^2
			&=-2\Real \int_{y_2}^R(\pt_{y_2}v)\ol v\df s\\
			&\leq 2\left(\int_{y_2}^Rs^\al|\pt_{y_2}v|^2\df s\right)^\f{1}{2}\left(\int_{y_2}^R s^{-\al} |v|^2\df s\right)^\f{1}{2}\leq 2y_2^{-\f{\al}{2}} \|v\|_{H_0^1((0,R);y_2^\al)}
		\end{split}
	\end{equation*}
	we deduce that 
	\begin{equation*}
		\begin{split}
			\int_I |v|^2\df y_2
			&\leq 2\|v\|_{H_0^1((0,R);y_2^\al)}\int_I y_2^{-\f{\al}{2}}\df y_2\\
			&\leq 2\|v\|_{H_0^1((0,R);y_2^\al)}\left(\f{1}{2}\mu^\f{1}{\ga}\right)^{-\f{\al}{2}}\f{1}{2}\mu^\f{\al}{2\ga}=2^\f{\al}{2}\|v\|_{H_0^1((0,R);y_2^\al)}. 
		\end{split}
	\end{equation*}
	If $I=\es$, its integral is zero and \eqref{11.20.10} alone suffices. 
	Combining these and \eqref{11.20.10} and \eqref{11.20.6}  we get 
	\begin{equation*}
		1\leq 4\|f\|_{L^2(0,R)}+\left[16\left(\f{3}{2}\right)^\f{\al}{2}+2^\f{\al}{2}\right]\|f\|_{L^2(0,R)}^\f{1}{2}. 
	\end{equation*}
	This estimate has the same form as \eqref{11.20.11}. These estimates also apply when $\mu^\f{1}{\ga}\geq R$. Thus, for all $\mu\in\R$ and for all $\ga\in \{1,2\}$, there exists a positive constant $C_\al\geq 1$, depending only on $\al$, such that 
	\begin{equation}\label{11.20.12}
		1\leq C_\al\left(\|f\|_{L^2(0,R)}+\|f\|_{L^2(0,R)}^\f{1}{2}\right). 
	\end{equation}

	If $\|f\|_{L^2(0,R)}\leq 1$, from \eqref{11.20.12}, we obtain 
	\begin{equation*}
		1\leq 2C_\al \|f\|_{L^2(0,R)}^\f{1}{2}\Longrightarrow \|f\|_{L^2(0,R)}\geq (2C_\al)^{-2}=:\ka. 
	\end{equation*}
If $\|f\|_{L^2(0,R)}>1$, then we have the same result. Combining these estimates with \eqref{11.20.5} and \eqref{11.20.7} yields
	\begin{equation}\label{11.20.13}
		\|(\mcB_R+\im \mu)v\|_{L^2(0,R)}\geq \f{1}{4C_\al^2} \|v\|_{L^2(0,R)}, \mbox{ for all } \mu\in \R \mbox{ and } v\in D(\mcB_R). 
	\end{equation}
	
	{\it Step 4}. From the imaginary-axis bound to decay.
	
	For $v_0\in D(\mcB_R)$ and $\tau>0$, set $v(t)=S_R(t)v_0$ and 
	\begin{equation*}
		z(t)=
		\begin{cases}
			\sin (\f{\pi t}{\tau})v(t), &0<t<\tau,\\
			0, &t\notin (0,\tau).  
		\end{cases}
	\end{equation*}
	Since $S_R(t)$ preserves $D(\mcB_R)$ and $\mcB_RS_R(t)v_0=S_R(t)\mcB_Rv_0$, we obtain 
	\begin{equation*}
		z\in H^1(\R; L^2(0,R))\cap L^2(\R; D(\mcB_R)) 
	\end{equation*}
	by extending  $z$ to be a function on $\R$ with $z(t)=0$ for all $t\notin (0,\tau)$. Then, in the sense of distributions, we obtain
	\begin{equation*}
		z_t+\mcB_Rz=g,\mbox{ with } g = \f{\pi}{\tau}\cos\left(\f{\pi t}{\tau}\right)v(t)\Ee_{(0,\tau)}(t), 
	\end{equation*}
where $\Ee_{(0,\tau)}$ is the characteristic function on $(0,\tau)$. Taking the Fourier transform in time gives $(\mcB_R+\im\zeta)\wh z(\zeta)=\wh g(\zeta)$. By \eqref{11.20.5} and the Fourier-Plancherel identity, we obtain
	\begin{equation*}
		\begin{split}
			\f{\pi^2}{\tau^2}\int_0^\tau \cos^2\left(\f{\pi t}{\tau}\right) \|v(t)\|_{L^2(0,R)}^2\df t\geq \ka^2\int_0 ^\tau \sin^2\left(\f{\pi t}{\tau}\right)\|v(t)\|_{L^2(0,R)}^2\df t. 
		\end{split}
	\end{equation*}
Since $S_R(t)$ is a contraction semigroup, $\|v(\tau)\|_{L^2(0,R)}\le\|v(t)\|_{L^2(0,R)} \le\|v_0\|_{L^2(0,R)}$ for $0<t<\tau$. Hence
	\begin{equation*}
		\ka^2\f{\tau}{2}\|S_R(\tau)v_0\|_{L^2(0,R)}^2\leq \f{\pi^2}{2\tau}\|v_0\|_{L^2(0,R)}^2. 
	\end{equation*}
	By the density of $D(\mcB_R)$ in $L^2(0,R)$, therefore,
	\begin{equation*}
		\|S_R(\tau)\|_{\mcL(L^2(0,R))}\leq \f{\pi}{\ka \tau}. 
	\end{equation*}
Choose $\tau_0=\f{2\pi}{\ka}$. Then $\|S_R(\tau_0)\|_{\mathcal L(L^2(0,R))}\le \frac 12$. For $t\ge0$, write $t=n\tau_0+\theta$, where $n=\lfloor \f{t}{\tau_0}\rfloor\in\{0,1,2,\ldots\}$ and $\theta\in[0,\tau_0)$. Here $\lfloor a\rfloor$ denotes the greatest integer less than or equal to $a$. The semigroup property yields
	\begin{equation}\label{11.21.1}
		\begin{split}
			\|S_R(t)\|_{\mcL(L^2(0,R))}\leq 2^{-n}\leq 2\exp\left(-\f{\ka \log 2}{2\pi}t\right), \mbox{ for all }t\geq 0. 
		\end{split}
	\end{equation}
	
	{\it Step 5}. Return to the original variables.
	
	Since $-\mcB_R$ generates the contraction semigroup $S_R(s)=e^{-s\mcB_R}$ on $L^2(0,R)$ and $U_\de: L^2(0,1) \to L^2(0,R)$ is a unitary operator satisfying $\mcA_k = |k|^\be U_\de^{-1} \mcB_R U_\de$ with $D(\mcA_k) = U_\de^{-1} D(\mcB_R)$ (by \eqref{11.20.1}), it follows from standard semigroup similarity and time-rescaling properties that $-\mcA_k$ generates a unique $C_0$-contraction semigroup $e^{-t\mcA_k}$ on $L^2(0,1)$ given by
	\begin{equation*}
		U_\de e^{-t\mcA_k} U_\de^{-1} = S_R(|k|^\be t),  \text{ for all } t \ge 0.
	\end{equation*}
Since $U_\de$ is unitary, we obtain
	\begin{equation*}
		\begin{split}
			\left\|e^{-t\mcA_k}u^0\right\|_{L^2(0,1)}\leq 2 \exp\left(-\f{\ka \log2}{2\pi} |k|^\be t\right) \|u^0\|_{L^2(0,1)}.  
		\end{split}
	\end{equation*}
	This proves \eqref{11.21.2} with $c_0=\f{\ka \log 2}{2\pi}$. The estimate is an operator norm bound and therefore holds for every $L^2$ initial datum, without any additional gradient assumption. The proof is complete.
\end{proof}

\begin{remark}
	When $\al=0$, the constant $\be$ in \eqref{11.05.2} is identical to that reported in \cite[Proposition 6.2]{Beauchard}.
\end{remark}

\subsection{Observability}

(11) Let $\wh\vp_k^0\in D(\mcA_k^2)$, $\ga=1$, $k\in\Z-\{0\}$, and $\wh u_k=0$. Taking $\lambda = \lambda_0 \geq 1$ in \eqref{11.04.6}, i.e., 
\begin{equation*}
	s = s_0 := (1+|k|^\f{1}{2})(1+T^8)\la_0 \geq 1,
\end{equation*}
from Proposition \ref{11.04.P1}, the bounds $\f{C_5}{T^8}\leq \Theta \leq \f{C_6}{T^8}$ on $(\f{T}{4},\f{3T}{4})$ (along with $\xi \leq C\Theta$), 
\begin{equation*}
	\Theta e^{-2s_0\Theta}\leq \Theta e^{-2\Theta}\leq C  \mbox{ for }\Theta\in (0,+\iy), 
\end{equation*}  
and Lemmas \ref{11.05.L1} and \ref{11.05.L2}, we obtain (noting $1+|k|^\f{1}{2}\leq 2|k|^\f{1}{2}$ for $k\in\Z-\{0\}$)
\begin{equation*}
	\begin{split}
		s_0\int_\La |\wh\vp_k(T)|^2\df x_2
		&\leq \f{C}{T}s_0e^{-\f{c_0}{2}|k|^\be T}\int_\f{T}{4}^\f{3T}{4}\int_{\La} |\wh\vp_k|^2\df x_2\df t\\
		&\leq \f{C}{T}s_0^2e^{-\f{c_0}{2}|k|^\be T}T^{16} e^{Cs_0\f{1}{T^8}} \int_\f{T}{4}^\f{3T}{4}\int_{\La}\Theta^2  |\wh\vp_k|^2e^{-2s_0\xi}\df x_2\df t\\
		&\leq CT^{15}e^{C|k|^\f{1}{2}\f{1+T^8}{T^8}-\f{c_0}{2}|k|^\be T}s_0^2\iint_{\La_T} \Theta^2 |\wh\vp_k|^2e^{-2s_0\xi}\df x_2\df t\\
		&\leq CT^{15}e^{C|k|^\f{1}{2}\f{1+T^8}{T^8}-\f{c_0}{2}|k|^\be T} s_0\int_0^T\Theta |\pt_{x_2}\wh\vp_k(1,t)|^2e^{-2s_0\Theta} \df t\\
		&\leq CT^{15}e^{C|k|^\f{1}{2}\f{1+T^8}{T^8}-\f{c_0}{2}|k|^\be T} s_0 \int_0^T|\pt_{x_2}\wh\vp_k(1,t)|^2\df t,  		
	\end{split}
\end{equation*}
where the positive constants $C, C_5, C_6$ depend only on $\al$. Dividing both sides by $s_0\geq 1$, we deduce that 
\begin{equation}\label{11.04.12}
	\int_\La |\wh\vp_k(T)|^2\df x_2\leq C_7T^{15}e^{C_8|k|^\f{1}{2}\f{1+T^8}{T^8}-\f{c_0}{2}|k|^\be T} \int_0^T|\pt_{x_2}\wh\vp_k(1,t)|^2\df t, 
\end{equation}
where the positive constants $C_7,C_8$ depend only on $\al$.

(12) Estimate \eqref{11.04.12} has so far been proved for
$\widehat\vp_k^0\in D(\mcA_k^2)$ only.
Its extension to arbitrary $L^2$ initial data will be obtained
after taking the frequency supremum, by applying the strong-data
estimate on a time interval separated from $t=0$.
We do not use an $L^2(\La)\to L^2(0,T)$ boundedness assertion
for the boundary derivative.

\begin{theorem}\label{11.07.T1}
	Let $T > 0$, $\al \in (0,1)$,  and $\ga = 1$  and  $\be = \frac{2-\al}{3-\al}$. Then, the homogeneous decoupled kinetic system 
	\begin{equation}\label{11.07.1}
		\begin{cases}
			\pt_t \wh\vp_k - \im k x_2^\ga \wh\vp_k - \pt_{x_2}\left(x_2^\al \pt_{x_2} \wh\vp_k\right) = 0, & \text{in } \La \times (0,T), \\
			\wh\vp_k(0,t) = 0, \  \wh\vp_k(1,t) = 0, & \text{for } t \in (0,T),
		\end{cases}
	\end{equation}
	has the following uniform boundary observability property.
	For each $k\ne0$ and $t>0$, its mild solution has a classical
	boundary derivative, by fixed-mode analytic smoothing.
	Define the nonnegative extended integral
	\[
	J_{k,T}(\widehat\vp_k^0)
	=\int_0^T|\partial_{x_2}\widehat\vp_k(1,t)|^2\,\mathrm{d}t
	\in[0,+\infty].
	\]
	There is $\mathcal C_{T,\alpha}>0$, independent of $k\ne0$,
	such that every mild solution with initial datum in $L^2(\La)$
	satisfies
	\begin{equation}\label{11.07.2}
		\|\widehat\vp_k(T)\|_2^2
		\le\mathcal C_{T,\alpha}J_{k,T}(\widehat\vp_k^0).
	\end{equation}
	This is a finite-valued observation inequality on
	$\{f\in L^2(\La):J_{k,T}(f)<\infty\}$ and an inequality of
	extended nonnegative quantities for arbitrary $L^2$ data.
\end{theorem}

\begin{proof}
	First consider strong initial data in $D(\mcA_k^2)$.
	For any interval length $L>0$, \eqref{11.04.12} gives
	\begin{equation}\label{11.07.3}
		\|\widehat\vp_k(L)\|_2^2
		\le C_7L^{15}e^{\mathcal H_L(|k|)}
		\int_0^L|\partial_{x_2}\widehat\vp_k(1,t)|^2
		\,\mathrm{d}t,
	\end{equation}
	where
	\[
	\mathcal H_L(K)=C_8K^{1/2}\frac{1+L^8}{L^8}
	-\frac{c_0}{2}K^\beta L.
	\]
	Because $\beta=(2-\alpha)/(3-\alpha)>1/2$, we have
	$\mathcal H_L(K)\to-\infty$ as $K\to\infty$.
	Thus
	\begin{equation}\label{11.07.4}
		\mathcal M_{L,\alpha}
		:=\sup_{k\in\mathbb Z\setminus\{0\}}\mathcal H_L(|k|)
		<\infty.
	\end{equation}
	Define $C^{\mathrm{sm}}_{L,\alpha}
	=C_7L^{15}e^{\mathcal M_{L,\alpha}}$.
	The strong-data observation constant is therefore uniform in $k$.
	
	Now let $\widehat\vp_k^0\in L^2(\La)$.
	For fixed $k$, $-\mcA_k$ is a bounded perturbation of the
	analytic generator $-\mcA_0$, so
	$e^{-\tau\mcA_k}\widehat\vp_k^0\in D(\mcA_k^2)$ for every
	$\tau>0$. Take $\tau=T/2$ and $L=T-\tau=T/2$.
	Apply the strong-data estimate to
	$v(t)=\widehat\vp_k(\tau+t)$ on $(0,L)$:
	\[
	\|\widehat\vp_k(T)\|_2^2
	\le C^{\mathrm{sm}}_{T/2,\alpha}
	\int_{T/2}^T|\partial_{x_2}\widehat\vp_k(1,t)|^2
	\,\mathrm{d}t
	\le C^{\mathrm{sm}}_{T/2,\alpha}J_{k,T}(\widehat\vp_k^0).
	\]
	The shifted observation integral is finite for each fixed mode.
	Taking $\mathcal C_{T,\alpha}=C^{\mathrm{sm}}_{T/2,\alpha}$
	proves \eqref{11.07.2}, with the stated interpretation if
	$J_{k,T}=+\infty$.
\end{proof}

\begin{theorem}\label{11.07.T2}
	Let $T > 0$, and $\al\in (0,1)$ and $\ga = 1$. Consider the global linear kinetic equation defined on the domain $\Omega = \mathbb{T}_{x_1} \times \La_{x_2}$:
	\begin{equation}\label{11.07.5}
		\begin{cases}
			\partial_t\vp-x_2^\gamma\partial_{x_1}\vp
			-\partial_{x_2}(x_2^\alpha\partial_{x_2}\vp)=0,
			&\text{in }\mathbb T\times\La\times(0,T),\\
			\vp(x_1,0,t)=\vp(x_1,1,t)=0,
			&\text{for }(x_1,t)\in\mathbb T\times(0,T).
		\end{cases}
	\end{equation}
	Then, the system \eqref{11.07.5} is globally observable from the boundary $x_2 = 1$. That is, for any initial condition $\vp^0 \in L^2(\mathbb{T} \times \La)$, the corresponding mild solution $\vp$ satisfies:
	\begin{equation*}
		\|\vp(\cdot, \cdot, T)\|_{L^2(\mathbb{T} \times \La)}^2 \leq \mathcal{M}(T, \al) \int_0^T \|\pt_{x_2}\vp(\cdot, 1, t)\|_{L^2(\mathbb{T})}^2 \df t,
	\end{equation*}
	where $\mathcal{M}(T, \al)$ is a positive constant depending exclusively on $T$ and $\al$, remaining uniform with respect to the solution trajectories. For general $L^2$ initial data, the observation integral is
	understood in the extended sense as
	\[
	\sum_{k\in\mathbb Z}\int_0^T
	|\partial_{x_2}\widehat\vp_k(1,t)|^2\,\mathrm dt
	\in[0,+\infty].
	\]
	When finite, the boundary trace is defined by the
	$L^2(\mathbb T\times(0,T))$ limit of its Fourier partial sums.
\end{theorem}

\begin{proof}
	We prove this theorem by the following steps. 
	
	{\it Step 1}. 
	We expand the solution as a Fourier series:
	\begin{equation}\label{11.07.6}
		\vp(x_1, x_2, t) = \sum_{k \in \Z} \wh\vp_k(x_2, t) \left(\f{1}{\sqrt{2\pi}}e^{\im k x_1}\right),
	\end{equation}
	where the Fourier mode $\wh\vp_k(x_2, t)$ represents the component associated with frequency $k$, see also \eqref{11.03.7}.
	
	{\it Step 2}. 
	By substituting the expansion \eqref{11.07.6} into the governing equation \eqref{11.07.5}, the orthogonality of the Fourier basis $\{\f{1}{\sqrt{2\pi}}e^{\im k x_1}\}_{k \in \Z}$ on $\mathbb{T}$ ensures that the system completely decouples into independent equations for each frequency mode.

	For $k = 0$, the zero-frequency component satisfies a purely dissipative parabolic equation without transport interaction:
	\begin{equation}\label{11.07.7}
		\begin{cases}
			\pt_t \wh\vp_0 - \pt_{x_2}\left(x_2^\al \pt_{x_2} \wh\vp_0\right) = 0, & \text{in } \La \times (0,T), \\
			\wh\vp_0(0,t) = 0, \ \wh\vp_0(1,t) = 0, & \text{for } t \in (0,T).
		\end{cases}
	\end{equation}
	For any non-zero frequency mode $k \in \Z - \{0\}$, the decoupled component satisfies exactly the framework of Theorem \ref{11.07.T1}:
	\begin{equation}\label{11.07.8}
		\begin{cases}
			\pt_t \wh\vp_k - \im k x_2^\ga \wh\vp_k - \pt_{x_2}\left(x_2^\al \pt_{x_2} \wh\vp_k\right) = 0, & \text{in } \La \times (0,T), \\
			\wh\vp_k(0,t) = 0, \ \wh\vp_k(1,t) = 0, & \text{for } t \in (0,T).
		\end{cases}
	\end{equation}

	{\it Step 3}. 
	First, we invoke Parseval's identity to express the total $L^2$-norm of the solution at the terminal time $T$ as the sum of the squares of the individual Fourier mode norms:
	\begin{equation}\label{11.07.9}
		\|\vp(\cdot, \cdot, T)\|_{L^2(\mathbb{T} \times \La)}^2 = \int_\mathbb{T} \int_\La |\vp(x_1, x_2, T)|^2 \df x_2 \df x_1 = \sum_{k \in \Z} \int_\La |\wh\vp_k(x_2, T)|^2 \df x_2.
	\end{equation}
	We separate the right-hand side of \eqref{11.07.9} into the zero-mode energy and the high-frequency summation:
	\begin{equation}\label{11.07.10}
		\|\vp(\cdot, \cdot, T)\|_{L^2(\mathbb{T} \times \La)}^2 = \int_\La |\wh\vp_0(x_2, T)|^2 \df x_2 + \sum_{k \in \Z - \{0\}} \int_\La |\wh\vp_k(x_2, T)|^2 \df x_2.
	\end{equation}
	
For the non-zero frequency components ($k \neq 0$), we directly apply the uniform quantitative observability inequality  \eqref{11.07.2} proved in Theorem \ref{11.07.T1}. This guarantees that for all $k \in \Z - \{0\}$, we have:
	\begin{equation}\label{11.07.11}
		\int_\La |\wh\vp_k(x_2, T)|^2 \df x_2 \leq \mathcal{C}_{T, \al} \int_0^T |\pt_{x_2}\wh\vp_k(1,t)|^2 \df t,
	\end{equation}
where the positive constant $\mathcal{C}_{T, \al} > 0$ depends exclusively on $T$ and $\al$, remaining strictly uniform and independent of the index $k$.

	{\it Step 4}. 
For the zero mode, first let $\widehat\vp_0^0\in D(\mcA_0^2)$. Apply Proposition~\ref{11.04.P1} with $k=0$, $\gamma=1$ and $s=(1+T^8)\lambda_0$. The weight $\Theta^2e^{-2s\xi}$ has a positive lower bound on $\overline\La\times[\f T4, \f{3T}{4}]$, while $\Theta e^{-2s\Theta}$ is bounded on $(0,T)$. Together with energy monotonicity, the Carleman estimate therefore yields zero-mode boundary observability. For arbitrary $L^2(\La)$ initial data, apply the same argument on $(\frac T2,T)$, using analytic smoothing at $\frac T2$, and enlarge the observation interval to $(0,T)$. Thus there exists $\mathcal C_0=\mathcal C_0(T,\alpha)>0$ such that
	\begin{equation}\label{11.07.12}
		\int_\La |\wh\vp_0(x_2, T)|^2 \df x_2 \leq \mathcal{C}_0 \int_0^T |\pt_{x_2}\wh\vp_0(1,t)|^2 \df t.
	\end{equation}
	
	Now, substituting both bounds \eqref{11.07.11} and \eqref{11.07.12} into the split energy relation \eqref{11.07.10}, we obtain:
	\begin{equation*}
		\begin{split}
			\|\vp(\cdot, \cdot, T)\|_{L^2(\mathbb{T} \times \La)}^2 
			&\leq \mathcal{C}_0 \int_0^T |\pt_{x_2}\wh\vp_0(1,t)|^2 \df t + \sum_{k \in \Z - \{0\}} \mathcal{C}_{T, \al} \int_0^T |\pt_{x_2}\wh\vp_k(1,t)|^2 \df t \\
			&\leq \max\{\mathcal{C}_0, \mathcal{C}_{T, \al}\} \left( \int_0^T |\pt_{x_2}\wh\vp_0(1,t)|^2 \df t + \sum_{k \in \Z - \{0\}} \int_0^T |\pt_{x_2}\wh\vp_k(1,t)|^2 \df t \right).
		\end{split}
	\end{equation*}
	We define the global physical observability constant as $\mathcal{M}(T, \al) \coloneqq \max\{\mathcal{C}_0, \mathcal{C}_{T, \al}\} > 0$, which depends solely on the fixed parameters $T$ and $\al$. Collapsing the summation inside the bracket via Tonelli's theorem yields:
	\begin{equation}\label{11.07.13}
		\|\vp(\cdot, \cdot, T)\|_{L^2(\mathbb{T} \times \La)}^2 \leq \mathcal{M}(T, \al) \int_0^T \left( \sum_{k \in \Z} |\pt_{x_2}\wh\vp_k(1,t)|^2 \right) \df t.
	\end{equation}

	If the observation is infinite, the conclusion is immediate. Otherwise, orthogonality gives convergence of the boundary Fourier series in $L^2(\mathbb T\times(0,T))$, and Parseval's identity yields, for almost every $t\in(0,T)$:
	\begin{equation}\label{11.07.14}
		\sum_{k \in \Z} |\pt_{x_2}\wh\vp_k(1,t)|^2 = \|\pt_{x_2}\vp(\cdot, 1, t)\|_{L^2(\mathbb{T})}^2.
		\end{equation}
			Substituting \eqref{11.07.14} directly into \eqref{11.07.13} leads to the desired global macro-observability framework:
			\begin{equation*}
				\|\vp(\cdot, \cdot, T)\|_{L^2(\mathbb{T} \times \La)}^2 \leq \mathcal{M}(T, \al) \int_0^T \|\pt_{x_2}\vp(\cdot, 1, t)\|_{L^2(\mathbb{T})}^2 \df t.
			\end{equation*}
			This completes the proof of Theorem \ref{11.07.T2}.
\end{proof}

\begin{remark}\label{11.07.R2}	
(1) For $\gamma=2$ and $\alpha\in(0,1)$, the exponent is
\begin{equation*}
	\beta = \frac{2-\alpha}{4-\alpha} < \frac{1}{2}
\end{equation*}
For each fixed $T>0$, the frequency profile produced by the present proof satisfies
\begin{equation*}
\mathcal H_T(K)
=C_8K^{1/2}\frac{1+T^8}{T^8}
-\frac{c_0}{2}K^\beta T\longrightarrow+\infty
\quad\text{as }K\to\infty.
\end{equation*}
Consequently, this argument does not yield an observability constant uniform in $k$ when $\gamma=2$. This divergence of the available upper bound does not prove that uniform boundary observability fails for that equation.	
	
(2) The results above concern observation on the whole upper
boundary $\mathbb T\times\{1\}$.
They do not establish internal observability from an arbitrary
nonempty open subset $\omega\subset\Omega$ for every $T>0$.
Because the horizontal transport speed is bounded and there is
no horizontal diffusion, such an extension requires additional
geometric and observation-time conditions and a separate proof.
Internal observability under suitable conditions is left open here.
\end{remark}

\section*{Conflict of interest}
All authors disclosed no relevant relationships.

\section*{Data availability}
No data was used for the research described in the article.

\end{document}